\documentclass{llncs}
\usepackage{amsmath, epsfig,graphicx}
\usepackage{amssymb,bbm}
\usepackage{latexsym,url}

\newcommand{\GG}{{\mathcal G}}
\newcommand{\DD}{{\mathcal D}}

\newcommand{\HH}{{\mathcal H}}

\newcommand{\oR}{{\mathbb R}}
\newcommand{\oN}{{\mathbb N}}
\newcommand{\EE}{\mathcal E}
\newcommand{\SSS}{{\mathcal S}}
\newcommand{\oZ}{{\mathbb Z}}

\newcommand{\bS}{\mathbf S}

\newcommand{\psd}{{\mathcal S}^n_+}

\newcommand{\PSD}{\mathcal{S}_{+}}
\newcommand{\EDM}{\text{\rm EDM}}
\newcommand{\fedm}{{\text{\rm ed}}}
\newcommand{\gd}{\text{\rm gd}}

\newcommand{\oG}{\overline{G}}

\newcommand{\rank}{{\text{\rm dim}}}
\newcommand{\ran}{{\text{\rm rank}}}
\newcommand{\tw}{{\rm tw}}

\newcommand{\qb}{\bold{q}}
\newcommand{\pb}{\bold{p}}
\newcommand{\p}{p}
\newcommand{\pp}{\bold{p}^{'}}
\newcommand{\q}{q}
\newcommand{\w}{w^{'}}
\newcommand{\y}{y}

\newcommand{\z}{z}

\newcommand{\Om}{\Omega}
\newcommand{\la}{\langle}
\newcommand{\ra}{\rangle}
\newcommand{\what}{\widehat}
\newcommand{\hG}{\widehat G}
\newcommand{\hE}{\widehat E}
\newcommand{\ignore}[1]{}
\newcommand{\Supp}{\SSS}

\begin{document}

\title{A new graph parameter related to  bounded rank positive semidefinite matrix completions}

\author{Monique Laurent\inst{1,2} \and Antonios Varvitsiotis\inst{1}}
\institute{Centrum Wiskunde \& Informatica (CWI), Amsterdam \and Tilburg University, The Netherlands.}

\maketitle

\begin{abstract}
The Gram dimension $\gd(G)$ of a graph $G$  is the smallest integer $k\ge 1$ such that any partial real symmetric matrix, whose entries are specified on the diagonal and at the off-diagonal positions corresponding to edges of $G$, can be completed to a positive semidefinite matrix of rank at most $k$ (assuming   a positive semidefinite completion exists).
For any fixed $k$ the class of graphs satisfying $\gd(G) \le k$ is minor closed, hence it can characterized by a finite list of forbidden minors.
We show that  the only minimal forbidden minor is $K_{k+1}$ for $k\le 3$ and  that there are two minimal forbidden minors: $K_5$ and $K_{2,2,2}$ for $k=4$.
We also show some close connections to Euclidean realizations  of graphs and to the graph parameter $\nu^=(G)$ of \cite{H03}.
In particular, our characterization of the graphs with $\gd(G)\le 4$   implies the forbidden minor characterization of the 3-realizable graphs of Belk and Connelly \cite{Belk,BC} and  of the graphs with $\nu^=(G) \le 4$  of van der Holst \cite{H03}. 
\end{abstract}


\section{Introduction}

Given a graph $G=(V=[n],E)$, a {\em $G$-partial matrix} is a real symmetrix $n\times n$ matrix whose entries are specified on the diagonal and at the off-diagonal positions corresponding to the edges of $G$. 
The problem of completing a partial  matrix  to a full positive semidefinite (psd) matrix is one of the most extensively studied  matrix completion problems.
A particular instance is   the completion problem for correlation matrices (where all diagonal entries are equal to 1) arising in probability and statistics, and it is also closely related to the completion problem for 
Euclidean distance matrices  with applications, e.g.,  to sensor network localization and  molecular conformation in chemistry.
We  give definitions below and refer, e.g.,  to \cite{DL97,L01} and further references therein for additional
 details.  

Among all   psd completions of a partial  matrix, the  ones with  the lowest possible rank are of particular importance. Indeed the rank of a matrix is often a good 
 measure of the complexity of the data it represents. 
As an  example, it is well known that the minimum  dimension of a Euclidean embedding of a  finite metric space can be expressed as the rank of an appropriate psd matrix (see e.g. \cite{DL97}).  Moreover, in applications, one is often interested in embeddings in low  dimension, say  2 or 3.
The problem of computing (approximate) low rank psd (or Euclidean) completions of a partial matrix  is a challenging non-continuous, non-convex problem which,
due to its great importance,   
has  been extensively studied (see, e.g., \cite{AAPW,AKW,RFP}, the recent survey \cite{KW} and further references therein).

\medskip
The following basic questions arise about psd matrix completions:  Decide whether a given partial rational matrix has a psd completion, what is the smallest rank of a completion, and if so find an (approximate) one (of smallest rank).
This leads to hard problems and of course the answer depends on the actual values of the entries of the partial matrix. 

However, taking a combinatorial approach to the problem and looking at   the structure of the graph $G$ of the specified entries, one can sometimes get tractable instances. For instance,  when the graph $G$ is chordal (i.e., has no induced circuit of length at least 4), the above questions are fully answered in \cite{GJSW84,La00} (see also the proof of Lemma \ref{lemcliquesum} below): There is a psd completion if and only if each fully specified principal submatrix is psd, the  minimum possible rank is equal to the largest rank of the fully specified principal submatrices, and such a psd completion can be found 
 in polynomial time (in the bit number model).  Further combinatorial characterizations (and some efficient algorithms for completions -- in the real number model) exist for graphs with no $K_4$-minor (more generallly when excluding certain splittings of wheels), see \cite{BJL96,La97,La00}.

\medskip In the present paper we focus on   the question of  existence of  low rank  psd completions. 
 Our approach is combinatorial, so we look for conditions on the graph $G$ of specified entries permitting to guarantee the existence of low rank completions. This is captured  by the  notion of {\em Gram dimension} of a  graph  which we introduce in Definition \ref{gramdef} below.
 
We use the following  notation:
 $\SSS^n$ denotes the set of symmetric $n\times n$ matrices and  $\SSS^n_+$ (resp., $\SSS^n_{++}$)
is the subset of all positive semidefinite (psd) (resp., positive definite) matrices.
For a matrix $X\in \SSS^n$, the notation $X\succeq 0$ means that  $X$ is psd.
Given  a graph  $G=(V=[n],E)$, 
 it will be convenient  to 
 identify $V$ with the set of diagonal pairs, i.e., to set     $V=\{(i,i)\mid i\in [n]\}$.
Then, a $G$-partial matrix corresponds to a vector $a\in \oR^{V\cup E}$ and $\pi_{V E}$ denotes the projection from $\SSS^n$ onto the subspace $\oR^{V\cup E}$  indexed by the diagonal entries and the edges of $G$.

\begin{definition}\label{gramdef}
The Gram dimension $\gd(G)$ of a graph $G=([n],E)$ is  the smallest integer $k\ge 1$ such that,
for any matrix $X\in \SSS^n_+$,
there exists  another matrix $X' \in \SSS^n_+$ with rank at most $k$ and such that
$\pi_{VE}(X)=\pi_{VE}(X^{'})$.
\end{definition}

Hence, if  a $G$-partial matrix  admits a psd completion,  it also has  one  of rank at most $\gd(G)$.
This motivates the study of bounds for  the graph parameter $\gd(G)$.
As we will see in Section~\ref{basic}, for any fixed $k$ the class of graphs with $\gd(G)\le k$ is closed under taking minors, hence it can be characterized  by a finite list of forbidden minors. Our main result is such a characterization for each integer $k\le 4$.\\ \vspace{0.10mm}

\noindent {\bf Main Theorem.} {\em For $k\le 3$,  $\gd(G)\le k$  if and only if  $G$ has no $K_{k+1}$ minor. For $k=4$,  $\gd(G)\le 4$ if and only if $G$ has no $K_5$ and $K_{2,2,2}$  minors. }\\ \vspace{0.10mm}

An equivalent way of rephrasing the notion of Gram dimension 
 is in terms of ranks of feasible solutions to semidefinite programs.
 Indeed, the  Gram dimension of a graph $G=(V,E)$ is at most $k$ if and only if the set
 $$S(G,a)=\{ X \succeq 0 \mid X_{ij}=a_{ij} \text{ for }  ij \in V \cup E\}$$
 contains a matrix of rank at most $k$ for all $a\in \oR^{V\cup E}$ for which $S(G,a)$ is not empty.
The set $S(G,a)$ is a typical instance of spectrahedron. Recall that a  
  {\em spectrahedron} is the convex region  defined as the intersection of the positive semidefinite cone with a finite set  of affine hyperplanes, i.e., the feasibility region 
  of a semidefinite program   in canonical  form:
\begin{equation}\label{sdpsform}\max \la A_0,X \ra \text{ subject to } \la A_j,X\ra=b_j,\   (j=1,\ldots,m), \qquad X \succeq 0.
\end{equation}
  If  the feasibility region
  of (\ref{sdpsform}) is not empty,
it follows from well known geometric results that it contains a 
matrix $X$ of rank $k$ satisfying 
${k+1\choose 2}\le m$ 
 (see e.g. \cite{Bar01}). 
 Applying this to the spectahedron $S(G,a)$,  
  we obtain the  bound
 $$\gd(G)\le \left\lfloor \frac{ \sqrt{1+8(|V|+|E|)}-1}{2}\right\rfloor.$$  For the complete graph $G=K_n$ the upper bound is equal to $n$, so  it is  tight. As we will see  one can get other  bounds depending on the structure of $G$; for instance, $\gd(G)$ is at most the tree-width plus 1 (cf. Lemma \ref{ktrees}). 

\medskip
As an application,  the Gram dimension can be used to bound the rank
 of optimal solutions to semidefinite programs. 
Namely,  consider  a semidefinite program 
in canonical form  (\ref{sdpsform}). Its {\em aggregated sparsity pattern} is the graph $G$  with node set $[n]$ and 
 whose edges are the pairs corresponding to the positions where at least one of the matrices  $A_j$ ($j\ge 0$) has a nonzero entry.
Then, whenever (\ref{sdpsform}) attains its maximum, it has  an optimal solution of rank at most $\gd(G)$. 
Results ensuring existence of low rank solutions are important, in particular, for approximation algorithms.
Indeed 
semidefinite programs are widely used as convex tractable relaxations to hard combinatorial problems. Then the rank one solutions typically correspond  to the  desired  optimal solutions of the discrete problem and low rank solutions can sometimes lead to improved performance guarantees (see, e.g., the result of \cite{AZ05} for max-cut and the result of \cite{BMZ02} for maximum stable sets).

As an illustration,
consider the max-cut problem for graph $G$ and its standard semidefinite programming relaxation:
\begin{equation}\label{maxcut}\max \frac{1}{4}\la L_G,X\ra \text{ subject to  }  X_{ii} =1\  (i=1,\ldots,n), \ \   X\succeq 0,
\end{equation}
where $L_G$ denotes the Laplacian matrix of $G$. Clearly, 
$G$ is the aggregated sparsity pattern of  the program 
(\ref{maxcut}). In particular, our main Theorem implies that if $G$ is $K_5$ and $K_{2,2,2}$ minor free, then (\ref{maxcut}) has an optimal
solution of rank at most four. (On the other hand recall 
that  the  max-cut problem  can be solved in polynomial time for $K_5$ minor free graphs  \cite{B83}).

In a similar flavor, for a graph $G=([n],E)$ with weights  $w\in \oR^{V\cup E}_+$,  the authors of \cite{GHW} 
study semidefinite programs of the form
$$\max \sum_{i=1}^n w_{i}X_{ii}\ \text{ s.t. } \sum_{i,j=1}^nw_{i} w_{j}X_{ij}=0,\
X_{ii}+X_{ij}-2X_{ij}\le w_{ij}\ (ij\in E),\ X\succeq 0,$$
and show
the existence of an  optimal solution of rank at most the tree-width of $G$ plus 1.
There is a large literature on dimensionality questions for various   geometric representations of graphs. We refer, e.g., 
 to 
 \cite{FH07,FH11,Hog08,Lo79,Lo01}  for results and further references. We will point out  links to the parameter $\nu^=(G)$ of \cite{H96,H03} in Section \ref{sec:Hein}.
 
\medskip
Yet another, more geometrical, way of interpreting the Gram dimension is in terms of isometric
 embeddings  in the spherical metric space~\cite{DL97}. 
 For this, consider the unit sphere $\bS^{k-1}=\{x \in \oR^k:  \|x\|=1\}$,  equipped with the distance 
$$d_{\bS}(x,y)=\arccos (x^Ty) \ \text{ for } x,y\in \bS^{k-1}.$$
Here, $\|x\|$ denotes the usual Euclidean norm.
Then $(\bS^{k-1},d_{\bS})$ is a metric space, known as  the {\em spherical metric space.}  
A graph $G=([n],E)$ has Gram dimension at most $k$ if and only if, for any assignment of vectors $p_1,\ldots,p_n \in \bS^h$ (for some $h\ge 1$),
there exists another assignment  $q_1,\ldots,q_n \in \bS^{k-1}$ such that 
$$d_{\bS}(p_i,p_j)=d_{\bS}(q_i,q_j), \text{ for } ij \in E.$$
In other words, this is the question of deciding whether a partial matrix can be realized in the $(k-1)$-dimensional spherical space.
The analogous question for the Euclidean metric space $(\oR^k,\|\cdot \|)$ has been extensively studied.  In Section~\ref{seclinks}  we will establish close connections with the notion of $k$-realizability of graphs introduced in~\cite{Belk,BC} and to the corresponding graph parameter $\fedm(G)$.

Complexity issues concerning the parameter $\gd(G,x)$ are discussed in~\cite{ELV12}. Specifically, given  a graph $G$ and a rational vector in $\EE(G)$, the  problem of deciding whether $\gd(G,x)\le k$ is proven to be NP-hard for every fixed $k\ge 2$~\cite{ELV12}.

\subsubsection*{Contents of the paper.} 
In Section~\ref{basic} we 
 determine   basic properties of the graph parameter $\gd(G)$ and in 
    Section~\ref{gramdim4} we reduce the proof of our main Theorem   to the problem of computing the Gram dimension of the two graphs $V_8$ and $C_5 \times C_2$. In Sections~\ref{seclinks} and \ref{sec:Hein}   we investigate the  links 
 of  $\gd(G)$ with the graph parameters  $\fedm(G)$ and  $\nu^=(G)$, respectively. Section \ref{sec:boundinggramdim} introduces the main ingredients for our proof: In Section~\ref{sec:genericity} we discuss some genericity assumptions we can make, in Section~\ref{sec:sdpformulation} we  show how to use semidefinite programming,
 in Section~\ref{ulemmas} we establish a number of useful lemmas, 
 and in Section~\ref{V8} we  show that $\gd(V_8)=4$.
Section~\ref{c2xc5} is dedicated to  proving  that $\gd(C_5\times C_2)=4$ -- this is the most technical part of the paper. Lastly, in Section~\ref{sec:concluding} we conclude with some comments and 
open problems.

\subsubsection*{Note.} Part of this work will appear as an extended abstract in the proceedings of ISCO 2012~\cite{LV12}.

\section{Preliminaries}
\subsection{Basic definitions and properties}\label{basic}
For a  graph $G=(V=[n],E)$ let  $ \PSD(G)=\pi_{VE}(\psd)\subseteq \oR^{V\cup E}$ denote the projection of the positive semidefinite cone onto $\oR^{V\cup E}$,  whose elements can be seen as the $G$-partial  matrices  
 that can be completed to a psd matrix. Let $\EE_n$ denote the set of matrices in $\SSS^n_+$ with an all-ones diagonal (aka the correlation matrices), and let
$\EE(G)=\pi_{E}(\EE_n)\subseteq \oR^E$ denote its projection onto the edge subspace $\oR^E$,  known as the {\em elliptope} of $G$;
we only project on the edge set since all diagonal entries are  implicitly known and equal to 1 for matrices in $\EE_n$.
\begin{definition}Given a graph $G=(V,E)$ and a vector $a\in \oR^{V \cup E}$,  a Gram representation of $a$ in $\oR^k$
consists of  a set of vectors $\p_1,\ldots,\p_n\in \oR^k$  such that $$\p_i^T\p_j=a_{ij}\  \forall  ij \in V \cup E.$$
The Gram dimension of   $a\in \PSD(G)$, denoted as $\gd(G,a)$,  is  the smallest integer $k$ for which $a$ has a Gram representation in $\oR^k$. 
\end{definition}

\begin{definition} 
The Gram dimension of a graph $G=(V,E)$ is defined as
\begin{equation}\label{gramdimdef}
 \gd(G)=\underset{a \in \PSD(G)}{\max} \gd(G,a).
 \end{equation}
\end{definition}
Clearly, the maximization in (\ref{gramdimdef}) can be restricted  to be taken over  all vectors $a \in \EE(G)$ (where all diagonal entries are implicitly taken to be equal to 1).
We  denote by  $\GG_k$  the class of graphs $G$  for which $\gd(G)\le k$. 

As a warm-up example, $\gd(K_n)=n$: The upper bound is clear   as $|V(K_n)|=n$ and the lower bound follows by considering,  e.g., $a=\pi_{V\cup E}(I_n)$.

\medskip We now investigate the behavior of the graph parameter $\gd(G)$ under some simple graph operations. Recall that   $G\backslash e$ (resp., $G\slash e$) denotes the graph  obtained from $G$ by deleting (resp., contracting) the edge $e$.
A graph $H$ is a {minor} of $G$ (denoted as $H\preceq G$) if $H$ can be  obtained from $G$ by successively deleting and contracting edges and deleting nodes.

\begin{lemma}\label{lemminor}
The graph parameter $\gd(G)$ is monotone nonincreasing with respect to edge deletion and contraction:
$\gd(G\backslash e), \gd (G\slash e)\le \gd(G) $ for any edge $e\in E$. 
\end{lemma}

\begin{proof} Let $G=([n],E)$ and $e\in E$. It follows directly from the definition that $\gd(G\backslash e)\le \gd(G)$. We show that 
$\gd(G\slash e)\le \gd(G)$. 
Say $e$ is the edge $(1,n)$ and $G\slash e=([n-1],E')$.
  Consider $ X \in \SSS_{+}^{n-1}$; we show that there exists $X'\in\SSS^{n-1}_+$ with rank at most $k=\gd(G)$ and such that $\pi_{E'}(X)=\pi_{E'}(X')$.
  For this, extend $X$ to the matrix $Y\in\SSS^n_+$ defined by $Y_{nn}=X_{11}$ and $Y_{in}=X_{1i}$ for $i\in[n-1]$.
  By assumption, there exists $Y'\in\SSS^n_+$ with rank at most $k$ such that $\pi_E(Y)=\pi_E(Y')$. Hence
  $Y'_{1i}=Y'_{ni}$ for all $i\in [n]$, so that 
  the principal submatrix $X'$ of $Y'$ indexed by $[n-1]$ has rank at most $k$ and satisfies
  $\pi_{E'}(X')=\pi_{E'}(X)$.
  \ignore{
   where $X=\left(\begin{array}{cc}
Y & b\\
b^T& 1
\end{array}\right)$. Define   $\what{X}=\left(\begin{array}{ccc}
 Y &  b& b\\
 b^T& 1& 1\\
 b^T & 1 & 1
 \end{array}\right)$ and notice that $\what{X} \in \psd$. Then there exists a matrix $X^{'} \in \psd$ with  $\ran X^{'} \le k$ such that $\pi_{E(G)}(\what{X})=\pi_{E(G)}(X^{'})$. Let $X^{'}[n]$ be the principal minor of $X^{'}$ defined by the first $n$ rows/columns. For $ij \in E(G/e)$ where $i,j\not=n-1$ we have that $X_{i,j}^{'}=X_{i,j}$.  Moreover if  $(n-1,i) \in  E(G/e)\cap E(G)$ then $X_{n-1,i}^{'}=\what{X}_{n-1,i}=X_{n-1,i}$. Lastly, for $(n-1,i) \in E(G/e)$ where $(n,i)\in E(G)$  we have that $X_{n-1,i}^{'}=X^{'}_{n,i}=\what{X}_{n,i}=X_{n-1,i}$.
 Summarizing, we have  that  $\pi_{E(G/e)}(X^{'}[n])=\pi_{E(G/e)}(X)$ and since $\ran X^{'}[n]\le k$ the claim follows. }
\qed\end{proof}

Let  $G_1=(V_1,E_1)$, $G_2=(V_2,E_2)$ be two graphs,  where  $V_1\cap V_2$  is a clique in both $G_1$ and $G_2$. Their {\em clique sum} is the graph $G=(V_1\cup V_2, E_1\cup E_2)$, also called their 
{\em clique $k$-sum} when $|V_1\cap V_2|=k$. The following result follows from well known arguments (cf. e.g. \cite{GJSW84};  a proof is included for completeness). For a matrix $X$ indexed by $V$ and a subset $U\subseteq V$, $X[U]$ denotes the principal submatrix of $X$ indexed by $U$.

\begin{lemma}\label{lemcliquesum}
If $G$ is the clique sum of two graphs $G_1$ and $G_2$,  then 
$$\gd(G)=\max\{\gd(G_1),\gd(G_2)\}.$$
\end{lemma}

\begin{proof}
The proof relies on the following fact: Two psd matrices $X_i$ indexed by $V_i$ ($i=1,2$) such that $X_1[V_1\cap V_2]=X_2[V_1\cap V_2]$ admit a common psd completion $X$ indexed by $V_1\cup V_2$ with rank $\max\{\rank (X_1),\rank (X_2)\}$. Indeed, let $u^{(i)}_j$ ($j\in V_i$) be a Gram representation of $X_i$ and let $U$ an orthogonal matrix mapping $u^{(1)}_j$ to $u^{(2)}_j$ for $j\in V_1\cap V_2$, then the Gram representation of $Uu^{(1)}_j$ ($j\in V_1$) together with $u^{(2)}_j$ ($j\in V_2\setminus V_1$) is such a common completion. \qed\end{proof}


Recall that  the {\em tree-width} of  a graph $G$, denoted by $\tw(G)$, is 
the minimum integer $k$ for which $G$ is contained (as a subgraph) in a clique sum of copies of ${K_{k+1}}.$  As a direct application of Lemmas~\ref{lemminor} and \ref{lemcliquesum} we obtain the following bound:

\begin{lemma}\label{ktrees} For any graph $G$, $\gd(G)\le \tw(G)+1$.

\end{lemma}



In view of Lemma \ref{lemminor}, the class $\GG_k$ of graphs with Gram dimension at most $k$ is closed under taking minors. 
Hence, by the celebrated graph minor theorem of \cite{RS}, 
it can be characterized by finitely many minimal forbidden minors. 

Clearly,  $K_n$ is a minimal forbidden minor for $\GG_{n-1}$ for all $n$, 
since contracting an edge yields   a graph with $n-1$ nodes and deleting an edge yields a graph with tree-width at most $n-2$.

It follows by its definition that  the tree-width of a graph  is a minor-monotone graph parameter.  
One can easily verify that 
$\tw(G) \le 1  \Longleftrightarrow   K_3 \not\preceq G$ and it is known that 
$\tw(G) \le 2   \Longleftrightarrow   K_4 \not\preceq G$ \cite{Du65}.
Combining these two facts with Lemma~\ref{ktrees} yields 
 the full list of forbidden minors for  the class $\GG_k$ when $k\le 3$.
\begin{theorem}\label{theosmallr}
For $k\le 3$, $\gd(G)\le k$  if and only if $G$ has no minor $K_{k+1}$.
\end{theorem}

\subsection{Characterizing graphs with Gram dimension  at most  4}\label{gramdim4}

The next natural question  is to characterize the class $\GG_4$.
Clearly, $K_5$ is a minimal forbidden minor for $\GG_4$.
We now  show that this is also the case for the complete tripartite graph $K_{2,2,2}$.
\ignore{
\begin{figure}[h]
\centering \includegraphics[scale=0.5]{figs/K222-ML}
\caption{The graph $K_{2,2,2}$.} 
\label{k222}
\end{figure}
}

\begin{lemma}\label{lemK222}
The graph $K_{2,2,2}$ is a minimal forbidden minor 
 for   $\GG_4$.
\end{lemma}

\begin{proof}
First we construct $a\in \EE(K_{2,2,2})$ with $\gd(K_{2,2,2},a)\ge 5$, thus implying $\gd(K_{2,2,2})\ge 5$. For this, let $K_{2,2,2}$ be obtained from $K_6$ by deleting the edges $(1,4)$, $(2,5)$ and $(3,6)$. Let $e_1,\ldots,e_5$ denote the standard unit vectors in $\oR^5$,
 let $X$ be the Gram matrix of the vectors 
$e_1,e_2,e_3,e_4,e_5$ and $(e_1+e_2)/\sqrt 2$ labeling the nodes $1,\ldots,6$, respectively, and let 
 $a\in \EE(K_{2,2,2})$ be the projection of $X$. We now verify that $X$ is the unique psd completion  of $a$ which shows  that $ \gd(K_{2,2,2},a)\ge 5$. Indeed the chosen Gram labeling of the matrix $X$ implies the following linear dependency:  $X[\cdot,6]=(X[\cdot,4]+X[\cdot,5])/\sqrt 2$ among its columns $X[\cdot,i]$ indexed respectively by $i=4,5,6$;
this implies that the unspecified entries $X_{14}, X_{25}, X_{36}$ are uniquely determined in terms of the specified entries of $X$.

On the other hand, one can easily verify that  $K_{2,2,2}$ is a partial 4-tree, 
therefore $\gd(K_{2,2,2})\le 5$. Moreover,  
deleting or contracting an edge in $K_{2,2,2}$ yields a partial 3-tree, thus with Gram dimension at most 4.
\ignore{
Next we verify that $K_{2,2,2}$ is a minimal forbidden minor for membership in $\GG_4$, i.e., 
  when deleting or contracting an edge one obtains a graph with Gram dimension at most 4.
Let $G$ be the graph obtained by deleting  the edge $(1,2)$ in $K_{2,2,2}$. If we add the edge $(3,6)$ to $G$, we obtain  
a graph which is the clique sum of three cliques on four nodes, namely the cliques on $\{3,4,5,6\}$, $\{2,3,4,6\}$, and on $\{1,3,5,6\}$. This implies that the Gram dimension of $G$ is at most 4.

Let $G$ be the graph obtained by contracting  the edge $(1,2)$ in $K_{2,2,2}$. As the subgraph of $G$ induced by $\{3,4,5,6\}$ is the clique sum of two $K_3$'s, its Gram dimension is at most 3 and thus the Gram dimension of $G$ is at most~4.
}		
\qed\end{proof}

By Lemma~\ref{ktrees}, all  graphs with tree-width at most three belong to   $\GG_4$. Moreover,  these graphs can be characterized in terms of forbidden minors as follows:


\begin{theorem}\cite{APC90} \label{theo3tree}
A graph $G$ has $\tw(G)\le 3$ if and only if $G$ does not have $K_5,K_{2,2,2}, V_8$ and $C_5 \times C_2$ as a minor.
\end{theorem} 

The graphs $V_8$ and $C_5\times C_2$ are shown   in Figures~\ref{v8} and~\ref{c2c5} below, respectively.
These four graphs   are   natural candidates for being forbidden minors for  the class $\GG_4$. We have already seen  that  for   $K_5$ and $K_{2,2,2}$ this is indeed the case. However, this is not true for  $V_8$ and $C_5 \times C_2$. Both   belong to $ \GG_4$, this will be proved in Section \ref{V8} for $V_8$ (Theorem \ref{theoV8})
and in Section \ref{c2xc5} for $C_5\times C_2$ (Theorem \ref{c2xc5thm}). These two results form the main technical part of the paper. Using them, we
can complete our characterization of the class $\GG_4$.

\begin{theorem}\label{theomain}
For a graph $G$, $\gd(G)\le 4$ if and only if $G$ does not have $K_5$ or $K_{2,2,2}$ as a minor.
\end{theorem}

\begin{proof}
Necessity follows from Lemmas \ref{lemminor} and \ref{lemK222}. Sufficiency
 follows from the following graph theoretical result, obtained by combining Theorem~\ref{theo3tree} with Seymour's splitter theorem (for a self-contained proof  see~\cite{H96}):
  every graph with no $K_5$ and $K_{2,2,2}$ minors can be obtained as a subgraph of a clique $k$-sum ($k\le 2$) of copies of graphs with tree-width at most 3, $V_8$ and $C_5 \times C_2$. Combining this fact with Theorems \ref{theoV8}, \ref{c2xc5thm} and Lemmas~\ref{lemcliquesum},  \ref{ktrees} the claim follows. \qed
\end{proof}


\subsection{Links to Euclidean  graph realizations}\label{seclinks}
In this section we investigate the links between  the Gram dimension and the notion of $k$-realizability of graphs introduced in~\cite{Belk,BC}. We start the discussion with some necessary definitions.  

Recall that a matrix $D=(d_{ij})\in \SSS^n$ is a {\em Euclidean distance matrix} (EDM) if there exist vectors $p_1,\ldots,p_n\in \oR^k$ (for some $k\ge 1$) such that $d_{ij}=\|p_i-p_j\|^2$ for all $i,j\in [n]$.
Then  $\EDM_n$ denotes the  cone of all $n \times n$ Euclidean distance matrices and, for a graph $G=([n],E)$, $
\EDM(G)=\pi_E(\EDM_n)$ is  the set of $G$-partial matrices that can be completed to a Euclidean distance matrix.

\begin{definition}\label{defedm}
Given a graph $G=([n],E)$ and  $d\in \oR_{+}^{E}$,  a Euclidean (distance) representation of $d$ in $\oR^k$
consists of  a set of vectors $\p_1,\ldots,\p_n\in \oR^k$  such that $$\| \p_i-\p_j\|^2=d_{ij}\  \forall  ij \in E.$$
Then, $\fedm(G,d)$ is the smallest  integer $k$ for which $d$ has a  Euclidean representation in $\oR^k$ and the graph parameter $\fedm(G)$ is defined as 
\begin{equation}\label{edmdimdef}
\fedm(G)=\underset{d \in \EDM(G)}{\max} \fedm(G,d).
\end{equation}
\end{definition}

In the terminology of~\cite{Belk,BC} a graph $G$ satisfying $\fedm(G)\le k$ is called {\em $k$-realizable}. 
It is easy to verify that the graph parameter $\fedm(G)$  is minor monotone. Hence  for any fixed $k\ge 1$ the class of graphs satisfying $\fedm(G)\le k$ can be characterized by a finite list of minimal forbidden minors. For $k\le 2$ the only forbidden minor is $K_{k+2}$. Belk and Connelly~\cite{Belk,BC} have determined the list of forbidden minors for $k=3$.
 
\begin{theorem}~\cite{Belk,BC} \label{theoBC}  
For a graph $G$, $\fedm(G)\le 3$ if and only if $G$ does not have $K_5$ and $K_{2,2,2}$  as minors. 
\end{theorem}

The hard part of the proof of~\cite{Belk,BC} is to prove  sufficiency, i.e., that if a graph $G$  has no $K_5$ and $K_{2,2,2}$ minors then $\fedm(G)\le 3$. We will obtain  this result  as a corollary of our main theorem (cf. Corollary~\ref{corol}). To this end, we  have to  establish some connections between the graphs parameters $\fedm(G)$ and $\gd(G)$.

There is a well known correspondence between psd and EDM completions (for  details and references see, e.g.,  \cite{DL97}).
Namely, for a graph  $G$,  let $\nabla G$ denote its {\em suspension graph},   obtained by adding a new node (the {\em apex} node, denoted by 0), adjacent to all nodes of $G$.
Consider the one-to-one  map $\phi:  \oR^{V \cup E(G)} \mapsto \oR_{+}^{E(\nabla G)}$,  which maps  $x\in \oR^{V \cup E(G)}$ to
 $d=\phi(x )\in \oR_{+}^{E(\nabla G)}$ defined  by
$$ d_{0i}= x_{ii}  \ (i\in [n]),\ \ \ d_{ij}=x_{ii}+x_{jj}-2x_{ij} \ (ij\in E(G)).$$
  Then the vectors $u_1,\ldots,u_n\in\oR^k$ form a Gram representation of $x$ if and only if the vectors  $u_0=0,u_1,\ldots,u_n$ form a Euclidean representation of $d=\phi(x)$ in $\oR^k$. This shows:

\begin{lemma}\label{covariance}
Let $G=(V,E)$ be a graph. Then,  $\gd(G,x)=\fedm(\nabla G,\phi(x))$ for any $x\in \oR^{V\cup E}$ and thus 
 $\gd(G)=\fedm(\nabla G)$.
 \end{lemma}


 For the Gram dimension of a graph one can show the following property:

 \begin{lemma}\label{lem1}
 Consider a graph $G=(V=[n],E)$ and its suspension graph   $\nabla G=([n]\cup\{0\},E\cup F)$, where $F=\{(0,i)\mid i\in [n]\}$. Given $x\in \oR^E$, its {\em $0$-extension} is the vector $y=(x,0)\in\oR^{E\cup F}$. If $x\in \SSS_+(G)$, then $y\in \SSS_+(\nabla G)$ and
 $\gd(\nabla G, y)=\gd(G,x)+1$. Moreover, 
 $\gd(\nabla G)= \gd(G)+1$.
 \end{lemma}

 \begin{proof}
 The first part is clear and implies $\gd(\nabla G)\ge \gd(G)+1$. 
 Set $k=\gd(G)$; we  show the reverse inequality $\gd(\nabla G)\le k+1$.
 For this, let $X\in \SSS^{n+1}_+$, written in block-form as $X=\left(\begin{matrix} \alpha & a^T \cr a & A\end{matrix}\right)$,
 where $A\in \psd$ and the first row/column is indexed by the apex node 0 of $\nabla G$.
 If $\alpha =0$ then $a=0$, $\pi_{VE}(A)$ has a Gram representation in $\oR^k$  and thus $\pi_{ V(\nabla G) E(\nabla G)}(X)$ too.
 Assume now $\alpha > 0$ and without loss of generality $\alpha =1$.
 Consider the Schur complement $Y$ of $X$ with respect to the entry $\alpha=1$, given by
  $Y=A-aa^T$.
  As $Y\in \SSS^n_+$,  there exists $Z\in \psd$ such that $\text{rank}(Z) \le k$ and $\pi_{V E}(Z)=\pi_{V E}(Y)$.
  Define the matrix $$X':= \left(\begin{matrix} 1 & a^T \cr a & aa^T\end{matrix}\right)+\left(\begin{matrix} 0 & 0 \cr 0 & Z\end{matrix}\right).$$
  Then, ${\rm rank}( X') ={\rm rank} ( Z)+1\le k+1$. Moreover, $X'$ and $X$ coincide at all diagonal entries as well as at all entries corresponding to edges of $\nabla G$. This concludes the proof that $\gd(\nabla G)\le k+1$.
  \qed\end{proof}

  We do not know whether the analogous property is true for the graph parameter $\fedm(G)$. On the other hand,  the following partial result holds, whose proof was communicated to us by A. Schrijver.

  \begin{theorem}\label{lemedmdim}
  For a graph $G$,
  $\fedm(\nabla G)\ge \fedm(G)+1$.
  \end{theorem}

  \begin{proof}
  Set $\fedm(\nabla G)=k$; we   show  $\fedm(G)\le k-1.$  We may assume that $G$ is connected (else deal with each connected component separately).
  Let  $d \in \EDM(G)$ and let $p_1=0,p_2, \ldots, p_n$ be a Euclidean representation of $d$ in $\oR^h$ ($h\ge 1$).
  Extend the $p_i$'s  to vectors $\widehat{p_i}=(p_i,0)\in \oR^{h+1}$  by appending an extra  coordinate equal to zero, 
   and set $\widehat{p}_0(t)=(0,t)\in \oR^{h+1}$ where $t$ is any positive real  scalar.
    Now consider  the distance   $\widehat{d}(t) \in \EDM(\nabla G)$ with Euclidean representation $\widehat{p_0}(t), \widehat{p_1},\ldots,\widehat{p_n}$.

       As $\fedm(\nabla G)=k$, there exists another Euclidean representation of $\widehat{d}(t)$  by  vectors $q_0(t), q_1(t),\ldots,q_n(t)$ lying  in  $\oR^k$.
         Without loss of generality, we can assume that
	   $q_0(t)=\widehat{p_0}(t)=(0,t)$ and $q_1(t)$ is the zero vector; 
	       for $i\in[n]$, write $q_i(t)=(u_i(t),a_i(t))$, where $u_i(t)\in \oR^{k-1}$ and $a_i(t) \in \oR$.
	           Then $\|q_i(t)\| = \|\widehat{p_i}\|=\|p_i\|$ whenever node $i$ is adjacent to node 1 in $G$.
		        As the graph $G$  is connected,  this implies that, for any $i\in [n]$,
			      the scalars   $\|q_i(t)\|$ ($t \in \oR_+$)  are bounded. Therefore    there exists a  sequence  $t_m \in \oR_+$ ($m\in\oN$)  converging to $+\infty$     and for which the sequence $(q_i(t_m))_m$ has a limit. Say  $q_i(t_m)=(a_i(t_m),u_i(t_m))$ converges to  $(u_i,a_i)\in \oR^k$ as $m \rightarrow +\infty$,  where $u_i\in\oR^{k-1}$ and $a_i\in\oR$.
				      The condition $\|q_0(t)-q_i(t)\|^2=\widehat{d}(t)_{0i}$ implies  that  $\|p_i\|^2+t^2=\|u_i(t)\|^2+(a_i(t)-t)^2$ and thus
				      $$  a_i(t_m)=\frac{a_i^2(t_m)+\|u_i(t_m)\|^2-\|p_i\|^2}{2t_m}  \hspace{0.2cm} \forall m\in \oN.$$
				       Taking the limit as $m \to \infty $ we obtain  that $\underset{m \to \infty }{\lim} a_i(t_m)=0$ and thus  $a_i=0$.
				        Then, for $i,j\in [n]$, $d_{ij}=\widehat{d}(t_m)_{ij}=\|(a_i(t_m),u_i(t_m))-(a_j(t_m),u_j(t_m))\|^2$ and taking the limit as $m \to +\infty$ we obtain that  $d_{ij}=\|u_i-u_j\|^2$.
					 This shows that  the vectors $u_1,\ldots,u_n$ form a Euclidean representation of $d$ in $\oR^{k-1}$.
					 \qed\end{proof}

Combining Lemma~\ref{covariance} with Theorem~\ref{lemedmdim} we obtain the following inequality relating the parameters $\fedm(G)$ and $\gd(G)$. 

\begin{theorem}\label{lem:ed_vs_gd} For any graph $G$ we have that  $\fedm(G)\le \gd(G)-1.$
\end{theorem}

Combining Theorem \ref{lem:ed_vs_gd} with our main theorem we can recover sufficiency in Theorem \ref{theoBC}.

\begin{corollary}\label{corol}For a graph $G$, if $G$ has no $K_5$ and $K_{2,2,2}$ minors then  $\fedm(G)\le 3$. 
\end{corollary}



\subsection{Relation with the graph parameter $\nu^=(G)$}\label{sec:Hein}
In this section we investigate the relation between  the Gram dimension of a graph  and  the graph parameter $\nu^=(G)$ introduced  in~\cite{H96,H03}.  
 Recall that the {\em corank} of a matrix $M\in \oR^{n \times n}$ is  the dimension of its kernel.  Consider the cone
 $${\mathcal C}(G)=\{M\in \SSS^n_+: M_{ij}=0\  \text{ for all distinct } i,j\in V \text{ with } (i,j)\not\in E\}$$
 which, as is well known, can be seen as the dual cone of the cone $\SSS_+(G)$.
We  now introduce the graph parameter $\nu^=(G)$.

\begin{definition}Given a graph $G=([n],E)$ the  parameter $\nu^=(G)$ is defined as the maximum corank of a matrix $M\in {\mathcal C}(G)$ satisfying the following  property: 
$$\forall X\in \SSS^n \ \ \ MX=0, \  X_{ii}=0 \   \forall i\in V, \ X_{ij}=0\  \forall (i,j)\in E \ \Longrightarrow X=0,$$
known as the {\em strong Arnold property}.
\end{definition}

It is proven in~\cite{H96,H03} that $\nu^=(G)$ is a minor monotone graph parameter. Hence  for any fixed integer $k \ge 1$ the class of graphs with $\nu^=(G)\le k$  
can be characterized by a finite family of minimal forbidden minors. For $k\le 3$ the only forbidden minor is $K_{k+1}$. Van der Holst~\cite{H96,H03} has determined the list of forbidden minors for $k=4$.

\begin{theorem}\label{thm:hein}\cite{H96,H03} For a graph $G$, $\nu^=(G)\le 4$ if and only if $G$ does not have $K_5$ and $K_{2,2,2}$ as minors.
\end{theorem}

By relating the two parameters $\gd(G)$ and $\nu^=(G)$ we can derive sufficiency in Theorem~\ref{thm:hein} from our main Theorem.
\begin{theorem}\label{thm:v_vs_gd} For any graph $G$, $\gd(G) \ge \nu^=(G)$.

\end{theorem}

\begin{proof} Let $k=\nu^=(G)$ be attained by some matrix $M\in \SSS^n_+$. Write $M=\sum_{i=1}^n \lambda_iv_iv_i^T$, where $\lambda_i\ge 0$, $\{v_1,\ldots,v_n\}$ is an orthonormal base of eigenvectors of $M$, and   $\{v_1,\ldots,v_k\}$ spans the kernel of $M$.
 Consider the matrix  $X=\sum_{i=1}^kv_iv_i^T$  and its projection $a=\pi_{E\cup V}(X)\in \SSS_+(G)$. By construction, ${\rm rank} (X)=k$. Hence   it is enough to show that $a$ has a unique psd completion, which will imply $\gd(G)\ge \gd(G,a)=k$.

For this let   $Y\in \SSS^n_+$ be another psd completion of $a$.
Hence   the matrix $X-Y$ has zero entries at all positions $(i,j)\in V\cup E$. 
Since the  matrix $M$ has zero entries at all off-diagonal positions  corresponding to non-edges of $G$,
 we deduce that  $\langle M,X-Y \rangle=0$. 
 On the other hand,  $\langle M,X\rangle=\sum_{i=1}^k\lambda_iv_i^TMv_i=0$. Therefore,   $\langle M,Y\rangle =0$. As  $M,X,Y$ are psd, the conditions
 $\langle M,X\rangle =\langle M,Y\rangle =0$ imply that $MX=MY=0$ and thus
 $M(X-Y)=0$.
Now we can apply the assumption that  the matrix $M$ satisfies  the strong Arnold property  and  deduce  that $X=Y$. \qed
\end{proof}

Combining Theorem~\ref{thm:v_vs_gd} with our main theorem we can recover  sufficiency in Theorem~\ref{thm:hein}. 

\begin{corollary}\label{cor2}For a graph $G$,  if $G$ does not have $K_5$ and $K_{2,2,2}$ as minors then $\nu^=(G)\le 4$.
\end{corollary}



Colin de Verdi\`ere \cite{dV98} studies  the graph parameter $\nu(G)$, defined as the maximum corank of a matrix $M$ satisfying the strong Arnold property and such that, for any  $i, j\in V$,  $M_{ij}=0 \Longleftrightarrow (i,j)\not\in E$. In particular he shows that $\nu(G)$ is unbounded for the  class of planar graphs. As $\nu(G)\le \nu^=(G)\le \gd(G)$, we obtain as a direct application:

\begin{corollary}The graph parameter $\gd(G)$ is unbounded for the class of planar graphs.

\end{corollary}

\section{Bounding the Gram dimension}\label{sec:boundinggramdim}

In this section we sketch our approach to show that $\gd(V_8)=\gd(C_5\times C_2)=4$.


\begin{definition} 
Given a graph $G=(V=[n],E)$, a   configuration of  $G$ is an assignment  of vectors $\p_1,\ldots,\p_n $  (in some space)  to the nodes of $G$;
 the pair  $(G,\pb)$ is called a  framework.
We use the notation   $\pb=\{\p_1,\ldots,\p_n\}$ and, for a subset $T\subseteq V$, $\pb_T=\{\p_i\mid i\in T\}$. Thus $\pb=\pb_V$ and we also set  
$\pb_{-i}=\pb_{V\setminus \{i\}}$.

 Two configurations  $\pb,\bold{q}$ of  $G$ (not necessarily lying  in the same space) are said to be  equivalent if $\p_i^T\p_j=q^T_i q_j$ for all $ij \in V\cup E$.

\end{definition}

Our objective  is to show that the two graphs $G=V_8$, $C_5 \times C_2$ belong to  $\GG_4$. That is, we must show that, given any $a\in \SSS_+(G)$, one can construct a Gram representation  $\qb$ of $(G,a)$ lying in the space $\oR^4$.

Along the lines of  \cite{Belk} (which deals with Euclidean distance realizations),  our strategy to achieve this is as follows:
First, we construct a `flat'  Gram representation $\pb$ of $(G,a)$ obtained by maximizing the inner product  $\p_{i_0}^T\p_{j_0}$ along a given pair $(i_0,j_0)$ which is not an edge of $G$. As suggested  in \cite{SY06}  (in the context of Euclidean distance realizations), this configuration  $\pb$ can be obtained by solving a semidefinite program; then 
  $\pb$ corresponds to  the Gram representation of an optimal solution $X$ to this program.
  
 In general we cannot yet claim that $\pb$ lies in $\oR^4$. However, we can derive useful information about $\pb$ by  using an optimal solution $\Om$ (which will correspond to a `stress matrix') to the dual semidefinite program. Indeed, the optimality condition $X\Om=0$ will imply some  linear dependencies among the $\p_i$'s  that can be used to show the existence of an equivalent representation $\qb$ of $(G,a)$ in low dimension.
 Roughly speaking, most often,   these dependencies  will force the majority  of the $\p_i$'s to lie in  $\oR^4$, and  one will be able to rotate each remaining vector $\p_j$   about the space  spanned by the vectors labeling  the neighbors of $j$ into $\oR^4$.  Showing that the initial representation $\pb$ can indeed be `folded' into $\oR^4$ as just described makes up the main body of the proof.

Before going into the details of the proof, we indicate  some  additional genericity assumptions that can be made w.l.o.g.  on the vector $a\in \SSS_+(G)$. This will be particularly useful when treating the graph $C_5\times C_2$.

\subsection{Genericity assumptions}\label{sec:genericity}

By definition, $\gd(G)$ is the maximum value of $\gd(G,a)$ taken over all $a\in \EE(G)$.
 Clearly we can restrict the maximum to be taken over all $a$ lying in a dense subset of $\EE(G)$.
 For instance, the set $\DD$ consisting of all $x\in \EE(G)$ that admit a positive definite completion in $\EE_n$ is dense in $\EE(G)$. We next identify a smaller dense subset $\DD^*$ of $\DD$ which will we use in our study of the Gram dimension of $C_5\times C_2$. 
 
 We start with a useful lemma, which characterizes the vectors $a\in \EE(C_n)$ admitting a Gram realization in $\oR^2$. Here $C_n$ denotes the cycle on $n$ nodes.
 
 \begin{lemma}\label{lemcycle}
Consider the vector $a=(\cos \vartheta_1, \cos \vartheta_2,\ldots,\cos\vartheta_n)\in \oR^{E(C_n)}$, where $\vartheta_1,\ldots,\vartheta_n\in [0,\pi]$. 
Then $\gd(C_n,a)\le 2$ 
 if and only if there exist
 $\epsilon \in \{\pm 1\}^n$ and $k\in \oZ$ such that 
$\sum_{i=1}^n\epsilon_i \vartheta_i=2k\pi$.
\end{lemma}

\begin{proof}
We prove the `only if' part. Assume that $u_1,\ldots,u_n\in \oR^2$ are unit vectors such that 
$u_i^Tu_{i+1}= \cos \vartheta_{i}$ for all $i\in [n]$ (setting $u_{n+1}=u_1$).
We may assume that  $u_1=(1,  0)^T$. 
Then, $u_1^Tu_2=\cos \vartheta_1$ implies that $u_2=(\cos (\epsilon_1 \vartheta_1),
\sin(\epsilon_1 \vartheta_1))^T$ for some $\epsilon_1\in\{\pm 1\}$.
Analogously, $u_2^Tu_3=\cos \vartheta_2$ implies  $u_3=
(\cos(\epsilon_1\vartheta_1+\epsilon_2\vartheta_2), \sin(\epsilon_1\vartheta_1+\epsilon_2\vartheta_2))^T$ for some $\epsilon_2\in \{\pm 1\}$.
Iterating, we find that there exists $\epsilon\in\{\pm 1\}^n$ such that 
$u_{i}=(\cos(\sum_{j=1}^{i-1}\epsilon_i \vartheta_i),  \sin (\sum_{j=1}^{i-1}\epsilon_i \vartheta_i))^T$ for $i=1,\ldots,n$.
Finally, the condition $u_n^Tu_1=\cos \vartheta_n= \cos (\sum_{i=1}^{n-1} \epsilon_i \vartheta_i)$ 
 implies $\sum_{i=1}^n\epsilon_i\vartheta_i\in 2\pi\oZ$.
 The arguments can be reversed to show the `if part'.
 \qed\end{proof}

 \begin{lemma} \label{lemD1}
Let $\DD^*$ be the set of all $a\in \EE(G)$ that admit a positive definite completion in $\EE_n$ 
satisfying the following condition:
 For any circuit $C$ in $G$, the restriction $a_C=(a_e)_{e\in C}$ of $a$ to $C$ does not admit a Gram representation in $\oR^2$.
 Then the set $\DD^*$ is  dense in $\EE(G)$.
 \end{lemma}

 \begin{proof}
 We show that $\DD^*$ is dense in $\DD$.
 Let $a\in \DD$ and set $a=\cos \vartheta$, where $\vartheta\in [0,\pi]^E$.
  Given a circuit  $C$  in $G$ (say of length $p$), it follows from Lemma~\ref{lemcycle}
  that  $a_C$ has a Gram realization in $\oR^2$ if and only if
  $\sum_{i=1}^p\epsilon_i\vartheta_i=2k\pi$ for some $\epsilon\in\{\pm 1\}^p$ and $k\in \oZ$ with $|k|\le p/2$.
  Let $\HH_C$ denote the union of the hyperplanes in $\oR^{E(C)}$ defined by these equations.
  Therefore,  $a\not\in \DD^*$
   if and only if $\vartheta \in \cup_C \HH_C$, where the union is taken over all circuits $C$ of $G$.
    Clearly we can find a sequence  $\vartheta^{(i)} \in [0,\pi]^E \setminus \cup_C \HH_C$ converging to $\vartheta$ as $i\rightarrow \infty$. 
    Then the sequence $a^{(i)}:=\cos \vartheta^{(i)}$   tends to $a$ as $i\rightarrow \infty$ and, for all $i$ large enough, $a^{(i)}\in \DD^*$.
    This shows that $\DD^*$ is a dense subset of $\DD$ and thus of $\EE(G)$.
    \qed\end{proof}

    \begin{corollary}\label{lemgeneric}
    For any graph $G=([n],E)$, $\gd(G)=\max \gd(G,a)$, where the maximum is over all $a\in\EE(G)$ admitting  a positive definite completion 
    and whose  restriction to any circuit of $G$ has no Gram representation in the plane.
    \end{corollary}

\subsection{Semidefinite programming formulation}\label{sec:sdpformulation}

We now describe how to model the  `flattening' procedure using semidefinite programming (sdp) and how to obtain a  `stress matrix' using  sdp duality.

Let $G=(V=[n],E)$ be a graph and let $e_0=(i_0,j_0)$ be a non-edge of $G$ (i.e., $i_0\ne j_0$ and $e_0\not\in E$).
Let $a\in \SSS_+(G)$ be a partial  positive semidefinite matrix for which we want  to show the existence of a Gram representation  in a small dimensional space.
For this consider the semidefinite program:
\begin{equation}\label{SDPP}
\max \ \langle E_{i_0j_0},X\rangle\ \ \text{\rm s.t. } \langle E_{ij},X\rangle =a_{ij} \ (ij\in V\cup E),\ 
\ X\succeq 0,
\end{equation}
where $E_{ij}=(e_ie_j^T+e_je_i^T)/2$ and  $e_1,\ldots,e_n$ are the standard unit vectors in $\oR^n$.
The     dual semidefinite program  of (\ref{SDPP})  reads:
\begin{equation}\label{SDPD}
\min  \sum_{ij\in V\cup E} w_{ij} a_{ij} 
\text{ \rm s.t. } \Omega= \sum_{ij\in V\cup E} w_{ij} E_{ij}-E_{i_0j_0}\succeq 0.
\end{equation}


\begin{theorem}\label{thmpsdstress} 

Consider a graph  $G=([n],E)$, a pair  $e_0=(i_0,j_0)\not\in E$, and let  $a\in \SSS_{++}(G)$.
 Then there exists a Gram realization $\pb=(p_1,\ldots,p_n)$ in $\oR^k$ (for some $k\ge 1$) of $(G,a)$ and a   matrix $\Omega=(w_{ij})\in \psd $  satisfying  
 \begin{equation}\label{patternOm0}
 w_{i_0j_0}\ne 0,
 \end{equation}
\begin{equation}\label{patternOm}
w_{ij}=0\ \text{ for all }  ij\not\in V\cup E\cup \{e_0\},
\end{equation}
\begin{equation}\label{eqstress1}
w_{ii}\p_i+\sum_{j|ij \in E \cup \{e_0\}} w_{ij}\p_j=0\ \text{for all } i\in [n],
\end{equation}
\begin{equation}\label{eqrk2}
\dim\la p_i,p_j\ra =2 \ \text{ for all } ij \in E.
\end{equation}
We refer to equation~(\ref{eqstress1}) as the {\em equilibrium condition}  at vertex $i$.
\end{theorem}

\begin{proof}
\ignore{
Let $E_{ij}$ denote the elementary symmetric matrices (with 1 at position $(i,j)$ and $(j,i)$ and 0 elsewhere). 
As described above, `flattening' corresponds to   searching  for an equivalent  configuration  which  maximizes the inner product on the non-edge $e_0$. This can be formulated as the following semidefinite program:
\begin{equation}\label{SDPP}
\max \langle E_{i_0j_0},X\rangle\ \ \text{\rm s.t. } \langle E_{ij},X\rangle =2x_{ij} \ (ij\in E),\ 
\langle E_{ii},X\rangle =x_{ii}\ (i\in [n]),\ X\succeq 0,
\end{equation}
The     dual semidefinite program  is:
\begin{equation}\label{SDPD}
\min \sum_{i=1}^nw_{ii}x_{ii} + 2\sum_{ij\in  E}w_{ij}x_{ij}\ \text{ \rm s.t. } \Omega= \sum_{ij\in E\cup\Delta} w_{ij} E_{ij}-E_{i_0j_0}\succeq 0.
\end{equation}
}
Consider the sdp (\ref{SDPP}) and its dual program (\ref{SDPD}). 
By assumption,  $a$ has a positive definite completion, hence  the  program  (\ref{SDPP}) is strictly feasible.
Clearly,  the dual program (\ref{SDPD}) is  also strictly feasible.
Hence there is no duality gap and the optimal values are attained in both programs.
  Let $(X,\Omega)$ be a pair of primal-dual optimal solutions.
Then  $(X,\Omega)$ satisfies the optimality condition:
$\langle X,\Omega \rangle =0$ or, equivalently,  $X\Om=0$.
Say $X$ has rank $k$ and let $\pb =\{\p_1,\ldots,\p_n\}\subseteq \oR^k$ be a Gram realization of $X$.
Now it suffices to observe that  the condition $X\Om=0$ can be reformulated as  the equilibrium conditions (\ref{eqstress1}). 
The conditions (\ref{patternOm0}) and (\ref{patternOm}) follow from the form of the dual program (\ref{SDPD}), and (\ref{eqrk2}) follows from the assumption $a\in \SSS_{++}(G)$.
\ignore{
Then, $\Omega^*X^*=0$ implies:
\begin{equation*}
w_{ii}\p_i+\sum_{j|ij \in E \cup e_0}^n w_{ij}\p_j=0,\ \text{for all } i\in [n].
\end{equation*}
Indeed, $\Omega^* X^*=0$ reads $\sum_j w_{ij}\p_j^T\p_k=0$ for all $i,k\in [v]$; hence, $\sum_jw_{ij}\p_j$ is orthogonal to all $\p_k$'s which implies that 
$\sum_jw_{ij}\p_j=0$ for all $i$.
}
\qed
\end{proof}


\ignore{
\begin{remark} \label{lemrank0}
Since we assume that $a\in \SSS_{++}(G)$,  the  configuration $\p$ returned  by Theorem~\ref{thmpsdstress} satisfies  the conditions: 
$\p_i\ne 0$ for all $i\in V$, and 
$\rank \la \p_i , \p_j\ra=2$  for any edge $(i, j) \in E$.
\end{remark}
}

Note that, using the following variant of  Farkas' lemma for semidefinite programming,  one can  show the existence of a nonzero positive semidefinite   matrix $\Om=(w_{ij})$ satisfying (\ref{patternOm}) and  the equilibrium conditions (\ref{eqstress1}) also   in the case when the sdp (\ref{SDPP}) is not strictly feasible,  however now with   $w_{i_0j_0}=0$. This remark will be useful in the exceptional case considered in Section \ref{secexceptionalcase} where we will have to solve again a semidefinite program of the form (\ref{SDPP}); 
however this program will have additional conditions imposing that some of the $p_i$'s are pinned so that one cannot anymore assume strict feasibility (see the proof of Lemma \ref{lemOmp}).

\begin{lemma} {\bf (Farkas' lemma for semidefinite programming)} \label{lemfarkas}  
(see   \cite{Lo95})
Let $b\in\oR^m$ and let $A_1,\ldots,A_m\in \SSS^n$ be given.
 Then exactly one of the following two assertions holds:
\begin{description}
\item[(i)] Either there exists  $X\in \SSS^n_{++}$ such that $\la A_j,X\ra =b_j$ for $j=1,\ldots,m$.
\item[(ii)] Or there exists a  vector $y\in \oR^m$ such that $\Om:=\sum_{j=1}^m y_j A_j\succeq 0$, $\Om\ne 0$  and $b^Ty \le 0$.
\end{description}
Moreover, for any  $X\succeq 0$ satisfying $\la A_j,X\ra =b_j$ ($j=1,\ldots,m$),   we have  in (ii) $\la X,\Om\ra = b^Ty=0$ 
and thus $X\Om=0$.\end{lemma}

\begin{proof}
Clearly,  if (i) holds then (ii) does not hold.
Conversely, assume (i) does not hold, i.e., $\SSS^n_{++}\cap {\mathcal L}=\emptyset$, 
where  $\mathcal L=\{X\in\SSS^n\mid \la A_j,X\ra=b_j\ \forall j\}$.
Then there exists a separating hyperplane, i.e., there exists a nonzero matrix $\Om\in\SSS^n$ and $\alpha\in \oR$ such that $\la\Om,X\ra\ge \alpha$ for all $X\in \SSS^n_{++}$ and $\la \Om,X\ra \le \alpha $ for all $X\in \mathcal L$.
This implies $\Om\succeq 0$, $\Om\in {\mathcal L}^\perp$, and $\alpha\le 0$, so that (ii) holds and the lemma follows.
\qed\end{proof}

\subsection{Useful lemmas}\label{ulemmas}

We start with some definitions about stressed frameworks and then we  establish some basic tools that we will repeatedly use later  in our proof for $V_8$ and $C_5\times C_2$. For a matrix $\Om\in\SSS^n$ its {\em support graph} is the graph $\SSS(\Om)$ is the graph with node set $[n]$ and with edges the pairs $(i,j)$ with $\Om_{ij}\ne 0$.

\ignore{
\begin{definition}{\bf (Support $\Supp(\Om)$  of a matrix $\Om$)}
Given a matrix $\Om=(w_{ij}) \in \SSS^n$, its {\em support graph} is the graph $\SSS(\Om)=(V_\Om,E_\Om)$, where $i\in V_\Om$ if $w_{ij}\ne 0$ for some $j\in [n]$ and  the edges are the pairs $ij$ with $w_{ij}\ne 0$. 
Clearly, if $\Om\succeq 0$, $i\in V_\Om$ if and only if $w_{ii}\ne 0$.
\end{definition}
}

\begin{definition}{ \bf (Stressed framework $(H,\pb,\Om)$) }
Consider a framework $(H=(V=[n],F), \pb)$.   A nonzero matrix $\Om=(w_{ij})\in \SSS^n$ is called a {\em stress matrix}  for the framework $(H,\pb) $ if  its support graph $\SSS(\Om)$  is contained in $H$ (i.e., $w_{ij}=0$ for all $ij\not\in V\cup F$) and 
 $\Om$ satisfies the equilibrium condition
 \begin{equation}\label{eqstress}
 w_{ii}\p_i+\sum_{j: (i,j) \in F}w_{ij} \p_j=0 \ \ \forall i\in V.
 \end{equation}
Then the triple  $(H,\pb,\Omega)$ is called a {\em stressed framework}, and a {\em psd stressed framework} if moreover $\Om\succeq 0$.

We let $V_\Om$ denote the set of nodes $i\in V$ for which $\Om_{ij}\ne 0$ for some $j\in V$. A node $i\in V$ is said to be a {\em 0-node} when $w_{ij}=0$ for all $j\in V$.  Hence, $V\setminus V_\Om$  is the set of all 0-nodes and, when $\Om\succeq 0$, $i$ is a 0-node if and only if $w_{ii}=0$.

The support graph $\SSS(\Om)$ of $\Om$ is called the {\em stressed graph}; its edges are called  the {\em stressed edges} of $H$ and the nodes $i\in V_\Om$ are called the  {\em stressed nodes}.

Given an integer $t\ge 1$, a node $i\in V$ is said to be a {\em $t$-node} if  its degree in the stressed graph $\SSS(\Om)$ is equal to $t$. 

\end{definition}

Throughout we will deal with stressed frameworks  $(H,\pb,\Om)$ obtained by applying Theorem \ref{thmpsdstress}. 
Hence the graph $H$ arises by adding a new edge $e_0$ to a given graph $G$, which we then denote as $H=\hG$, as indicated below.

\begin{definition}
{ \bf (Extended graph $\hG$)}
Given a graph $G=(V=[n],E)$ and a fixed pair  $e_0=(i_0,j_0)$ not belonging to $E$,
 we set $\hG=(V,\hE=E\cup\{e_0\})$.
\end{definition}

We now group some useful lemmas which we will use in order to show that a given framework $(H,\pb)$ admits an equivalent 
configuration  in  lower dimension.

Clearly, the stress matrix  provides some linear dependencies among the vectors $\p_i$  labeling the stressed nodes, but it gives no information about the vectors  labeling the 0-nodes. 
However, if we have a set $S$ of 0-nodes forming a stable set, then we can use the following lemma in order to `fold'  the corresponding vectors $\p_i$ ($i\in S$) in a lower dimensional space.  

\begin{lemma}\label{folding} {\bf (Folding the vectors labeling a stable set})
Let $(H=(V,F),\pb)$ be a  framework and let $T\subseteq V$. Assume that  $S=V\setminus T$ is a stable set in $H$, 
that each node $i\in S$ has degree at most $k-1$ in $H$, and that  $\rank\la \pb_T\ra\le k$.
Then there exists a configuration $\bold{q}$ of $H$  in $ \oR^k$ which is equivalent to $(H,\pb)$.
\end{lemma}

\begin{proof}
Fix a node $i\in S$. Let $N[i]$ denote the closed neighborhood of $i$ in $H$ consisting of $i$ and the nodes adjacent to $i$.
By assumption, $|N[i]| \le k$ and   both sets of vectors $\pb_T$ and  $\pb_{N[i]}$ have rank at most $k$. Hence  one can find an orthogonal matrix $P$
 mapping all vectors $p_j$ ($j\in T\cup N[i]$) 
  into the space $\oR^k$. 
Repeat this construction  with every other node of $S$.  As no two nodes of $S$ are adjacent, this produces a configuration $\qb$ in $\oR^k$ which is equivalent to $(H,\pb)$.
\qed

\end{proof}


The next lemma uses the stress matrix to upper bound  the dimension of a given stressed configuration.

\begin{lemma}\label{base} { \bf (Bounding the dimension)} 
 Let $(H=(V=[n],F),\pb,\Omega)$ be a psd stressed framework. 
Then $\rank \la \pb_V\ra \le n-2$, except $\rank \la \pb_V\ra \le n-1$ if $\SSS(\Om)$ is a clique.
\end{lemma}
\ignore{
\begin{itemize}
\item[(i)] $\dim\langle\p_1,\ldots,\p_n\rangle\le n-1$.
\item[(ii)] If $G\not=K_n$ then $\dim\langle\p_1,\ldots,\p_n\rangle \le n-2.$
\end{itemize}
}

\begin{proof}
Let $X$ denote the Gram matrix of the $\p_i$'s, so that $\text{rank} (X)=\rank \la \pb_V\ra$. 
By assumption,  $X\Om=0$. This implies that $\text{rank} (X)\le n-1$.  Moreover, if $\SSS(\Om)$ is not a clique, then $\text{rank} (\Om)\ge 2$ and thus $\text{rank} (X)\le n-2$.
\qed \end{proof}


The next lemma indicates how  1-nodes can occur  in a stressed framework.

\begin{lemma}\label{lem1node}
Let $(H=(V,F),\pb,\Om)$ be a stressed framework.
  If node $i$ is a 1-node in the stressed graph $\SSS(\Om)$, i.e., there is a unique edge $ij\in F$ such that $w_{ij}\ne 0$, then 
$\rank\la\p_i, \p_j\ra \le 1$. 
\end{lemma}

\begin{proof}
Directly, using the equilibrium condition (\ref{eqstress}) at node $i$. 
\qed\end{proof}
We now consider   2-nodes in a  stressed framework.
First recall the notion of Schur complement.  For a matrix $\Om=(w_{ij})\in \SSS^n$ and $i\in [n]$ with $w_{ii}\ne 0$, the {\em Schur complement} of $\Om$ with respect to its $(i,i)$-entry is the matrix, denoted as $\Om_{-i}=(w'_{jk})_{j,k\in [n]\setminus \{i\}}\in \SSS^{n-1}$, with entries
$w'_{jk}= w_{jk} -w_{ik}w_{jk}/w_{ii}$ for $i,j\in [n]\setminus \{i\}$.
As is well known, $\Om\succeq 0$ if and only if $w_{ii}>0$ and $\Om_{-i}\succeq 0$.
We also need the following notion of `contracting a degree two node' in a graph.

\begin{definition}
Let $H=(V,F)$ be a graph, let   $i\in V$ be a node of degree two in $H$ which is adjacent to nodes $i_1,i_2\in V$.
The graph obtained by {\em contracting node $i$} in $H$ is the graph $H/i$ with node set 
$V\setminus \{i\}$ and with edge set
$F/i= F\setminus \{(i,i_1),(i,i_2)\} \cup \{(i_1,i_2)\}$ (ignoring multiple edges).
\end{definition}


\begin{lemma}\label{contract2node} {\bf (Contracting a 2-node)}
Let $(H=(V,F),\pb, \Omega)$ be a psd stressed framework, 
let   $i\in V$ be a 2-node in the stressed graph $\Supp(\Om)$ and set $N(i)=\{i_1,i_2\}$.
Then $\p_i\in \la \p_{i_1},\p_{i_2}\ra$ and thus $\rank \la \pb\ra = \rank \la \pb_{-i}\ra.$

Moreover, if  the stressed graph $\Supp(\Om)$ is not  the complete graph  on  $N[i]=\{i,i_1,i_2\}$, then
$(H/i,\pb_{-i},\Omega_{-i})$ is a psd stressed framework. 
\end{lemma}

\begin{proof}
 The equilibrium condition at node  $i$ implies $\p_i \in \la \p_{i_1},\p_{i_2}\ra$.
Note that the Schur complement  $\Om_{-i}$ of $\Omega$ with respect to the $(i,i)$-entry $w_{ii}$ has entries 
$w'_{i_1i_2}=w_{i_1i_2}- w_{ii_1}w_{ii_2}/w_{ii}$, 
$w_{i_ri_r}'=w_{i_ri_r}- w_{ii_r}^2/w_{ii}$ for $r=1,2$, and $w'_{jk}=w_{jk}$ for all other edges $jk$ of $H/i$.
As $\Om\succeq 0$ we also have   $\Om_{-i}\succeq 0$. Moreover, $\Om_{-i}\ne 0$. Indeed, $w_{i_1i_2}'\ne 0$ if $(i_1,i_2)\not\in F$; 
otherwise, as $\SSS(\Om)$ is not the clique on $N[i]$, 
there is another edge $jk$ of $H/i$ in the support of $\Om$ so that $w_{jk}'=w_{jk}\ne 0$.

In order  to show that $\Om_{-i}$ is a stress matrix for $(H/i,\pb_{-i})$,
 it suffices to check the stress equilibrium at the nodes $i_1$ and  $i_2$.  To fix ideas consider node  $i_1$. Then  we can rewrite 
$w'_{i_1i_1}\p_{i_1}+w'_{i_1i_2}\p_{i_2}+\sum_{j\in N(i_1)\setminus \{i_2\}}w'_{i_1j}\p_j$ as  
$
( \sum_{j} w_{i_1j}\p_j) - \left(w_{ii}\p_i+ w_{ii_1}\p_{i_1}+w_{ii_2}\p_{i_2}\right) w_{ii_1}/w_{ii},$
where both terms are equal to 0 using the equilibrium conditions of $(\Om,\pb)$ at nodes $i_1$ and $i$.
\ignore{
Finally we check the last statement about 0- and 1-nodes. Clearly, if $k\ne i_1,i_2$, then  $k$ has the same degree in $\SSS(\Om)$ as in $\SSS(\Om_{-i})$. Consider now (say) node $k=i_1$.
If the pair $(i_1,i_2)$ is  not stressed in $\SSS(\Om)$  then  node  $i_1$ has again the same degree in $\SSS(\Om)$ and in $\SSS(\Om_{-i})$.  Otherwise, $i_1$ has degree $d-1$ in $\SSS(\Om_{-i})$ if its degree is $d$   in $\SSS(\Om)$. 
Then, $d=1$ implies $\rank \la \p_{i_1},\p_i\ra\le 1$, while $d=2$ implies $(i_1,i_2)\in F$ and $\rank \la \p_{i_1},\p_{i_2}\ra\le 1.$
}
\qed
\end{proof}

We will apply the above lemma iteratively to contract a set $I$ containing  several  2-nodes. 
Of course, in order to obtain useful information, we want to be able to claim that, after contraction, we obtain a stressed framework $(H/I, \pb_{V\setminus I}, \Om_{-I})$, i.e., with $\Om_{-I}\ne 0$.
Problems might occur if at some step we get a stressed graph which is a clique on 3 nodes. Note that this can happen only when a connected component of the stressed graph is a circuit.  However, when we will apply this operation of contracting 2-nodes to the case of $G=C_5\times C_2$, we will make sure that this situation cannot happen; that is, we will show that we may assume  that the stressed graph does not have a connected component which is a circuit (see  Remark \ref{remcontract2node} in Section~\ref{uaddlemmas}).

\subsection{The graph $V_8$ has Gram dimension 4}\label{V8}
\begin{figure}[h] 
\centering \includegraphics[scale=0.5]{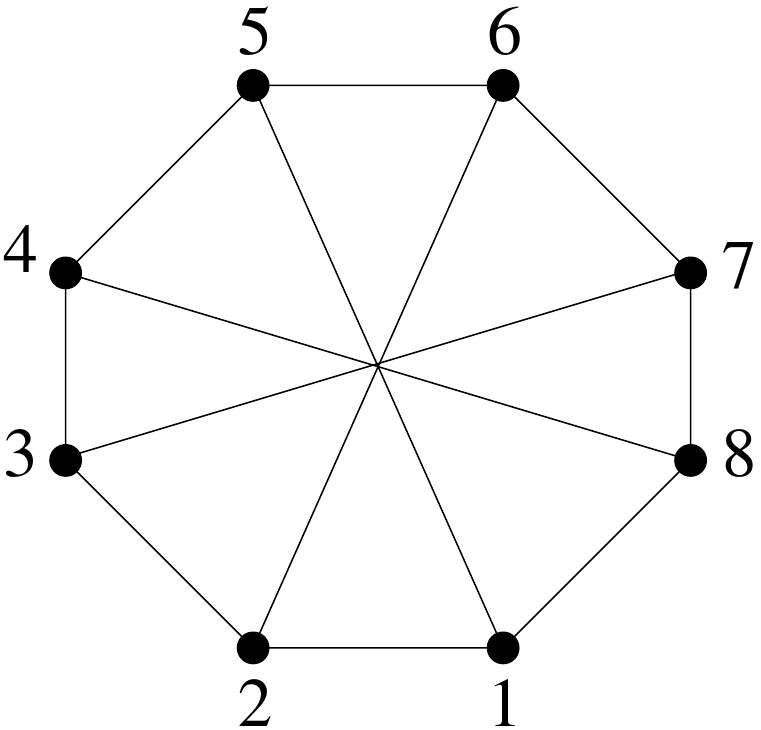}
\caption{The graph $V_8$.}
\label{v8} 
\end{figure}

Let $V_8=(V=[8],E)$ be the graph  shown in Figure~\ref{v8}. 
In this section we use the tools developed above to show that  $V_8$ has Gram dimension 4.

\begin{theorem}\label{theoV8}
The graph $V_8$ has Gram dimension 4.
\end{theorem}

\begin{proof} Set $G=V_8=([8],E)$. Clearly $\gd(G)\ge 4$ since $K_4$ is a minor of $ G$.
Fix $a\in \SSS_{++}(G)$; we show that $(G,a)$ has a Gram realization in $\oR^4$.
For this we first apply Theorem \ref{thmpsdstress}. 
As stretched edge $e_0$, we choose the pair $e_0=(1,4)$ and we denote by $\hG=([8], \hE=E\cup\{(1,4)\})$ the extended graph obtained by adding the stretched pair  $(1,4)$ to $G$.
 Let $\pb$ be the initial Gram realization of $(G,a)$ and let $\Om=(w_{ij})$ be the corresponding stress matrix obtained by applying Theorem \ref{thmpsdstress}.
We now show how  to construct from $\pb$  an equivalent realization $\q$ of $(G,a)$ lying in $\oR^4$.

In view of  Lemma \ref{folding}, we know that we  are done if we can find a subset $S\subseteq V$ which is stable  in the graph $G$ and satisfies  $\rank \la \pb_{V\setminus S}\ra\le 4$. This permits to deal with 1-nodes. Indeed suppose  that there is a 1-node in the stressed graph $\SSS(\Om)$. In view of Lemma \ref{lem1node} and~(\ref{eqrk2}), this can only be node 1 (or node 4) (i.e., the end points of the stretched pair) and  $\rank \la \p_1,\p_4\ra \le 1$. Then, choosing the stable set  $S=\{2,5,7\}$, we have  $\rank\la \pb_{V\setminus S}\ra \le 4$ and we can conclude using Lemma \ref{folding}.
 Hence we can now assume that there is no 1-node in the stressed graph $\SSS(\Om)$.

\medskip
Next, observe that we are done in any of the following two cases:
\begin{itemize}
\item [(i)]
There exists a set $T\subseteq V$ with $|T|=4$ and $\rank \la \pb_T\ra \le 2$.
\item[(ii)]
There exists a set $T\subseteq V$ of cardinality $|T|=3$ such that  $T$ does not consist of three consecutive nodes on the circuit $(1,2,\ldots,8)$ and $\rank \la \pb_T\ra \le 2$.
\end{itemize}
Indeed, in case (i) (resp., case (ii)),   there is a stable set $S \subseteq V\setminus T$ of cardinality $|S|=2$ (resp., $|S|=3$), so that $|V\setminus (S\cup T)|=2$ and thus 
$\rank \la \pb_{V\setminus S} \ra \le \rank \la \pb_T\ra + \rank \la \pb_{V\setminus (S\cup T)}   \ra \le 2+2 =4$.

\medskip
Hence we may assume that we are not in the situation of cases (i) and (ii).

\medskip
Assume first that one of the nodes in $\{5,6,7,8\}$ is a 0-node. Then all of them are 0-nodes. Indeed, if (say) 5 is a 0-node and 6 is not a 0-node then  the equilibrium equation at node 6 implies that $\rank \la \p_6,\p_7,\p_2\ra \le 2$ so that we are in the situation of case (ii). As nodes 1, 4 are not 1-nodes,  the stressed graph $\SSS(\Om)$  is the circuit $(1,2,3,4)$. Using Lemma \ref{contract2node}, we deduce that   $\rank \la \p_1,\p_2,\p_3,\p_4\ra \le 2$ and thus we are in the situation of case (i) above.

\medskip
Assume now that none of the nodes in $\{5,6,7,8\}$ is a 0-node but one of the nodes in $\{2,3\}$ is a 0-node. Then both nodes 2 and 3 are 0-nodes (else we are in the situation of case (ii)).
Therefore, both  nodes $6$ and $7$ are 2-nodes. Applying Lemma \ref{contract2node}, after contracting both nodes 6,7, we obtain a stressed framework  on $\{1,4,5,8\}$ and thus 
$\rank\la  \pb_{V\setminus\{2,3\}}\ra = \rank\la \p_1,\p_4,\p_5,\p_8 \ra$.
 Using Lemma \ref{base},  we deduce that 
$\rank \la \p_1,\p_4,\p_5,\p_8\ra \le 3$.
Therefore, $\rank \la \pb_{V\setminus \{3\}}\ra \le 4$ and one can find a new realization $\qb$ in $\oR^4$ equivalent to $(G,\pb)$ using  Lemma \ref{folding}.

\medskip
Finally assume that  none of the nodes in $\{2,3,5,6,7,8\}$ is a 0-node.
We show that  $\la {\bf p}\ra = \la \p_2,\p_3,\p_6,\p_7\ra$. 
The  equilibrium equation at node 6 implies that $\rank \la \p_2,\p_5,\p_6,\p_7 \ra \le 3$. Moreover, $\rank \la \p_2,\p_6,\p_7\ra=3$ (else we are in case (ii) above). Hence $\p_5\in \la \p_2, \p_6,\p_7\ra$.
Analogously, the equilibrium equations at nodes 7,2,3 give that $\p_8,\p_1,\p_4 \in  \la \p_2,\p_3, \p_6,\p_7\ra$, respectively.
\qed
\end{proof}

\section{The graph $C_5\times C_2$ has Gram dimension 4}\label{c2xc5}

\begin{figure}[h] 
\centering \includegraphics[scale=0.5]{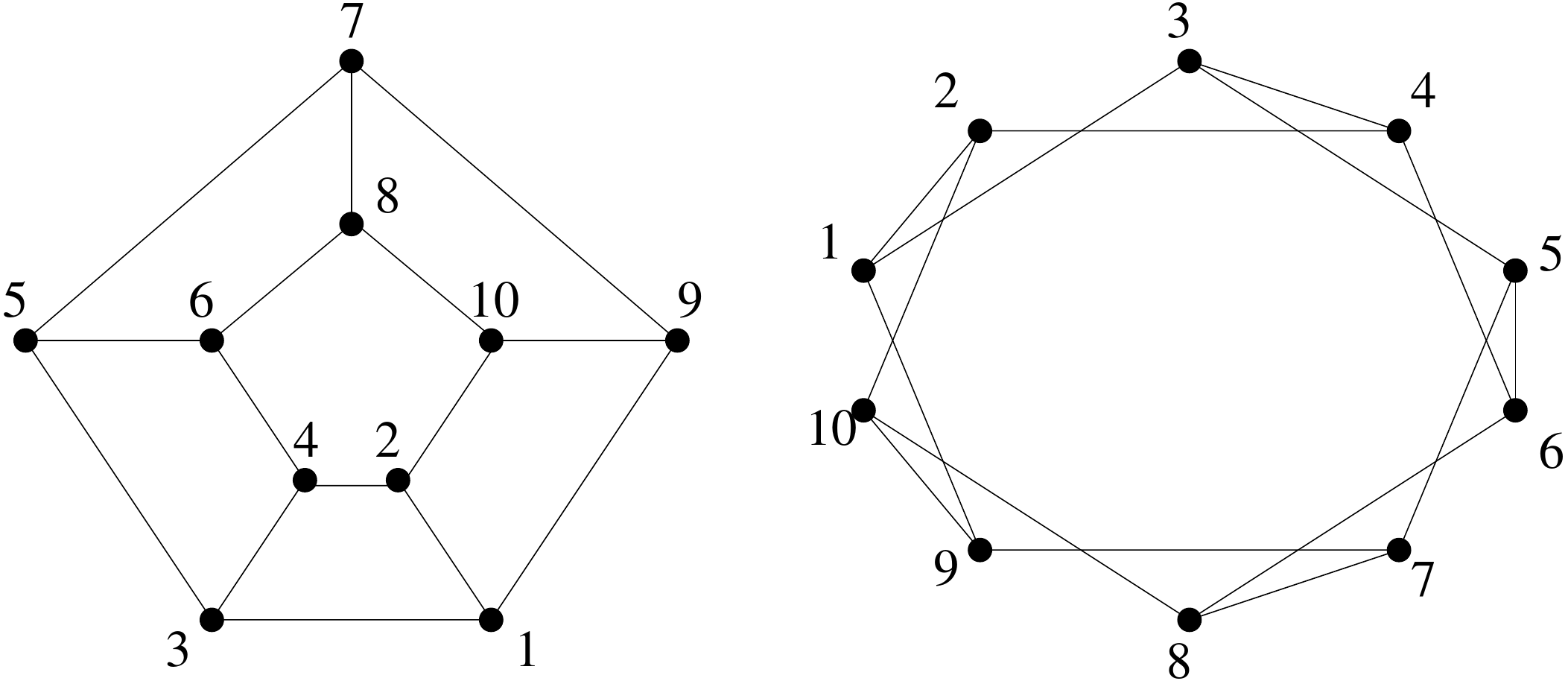}
\caption{Two drawings of the graph $C_5\times C_2$.}
\label{c2c5} 
\end{figure}

This section is devoted to proving that the graph $C_5\times C_2$ 
has Gram dimension 4. The analysis  is considerably more involved than the analysis for $V_8$.
Figure~\ref{c2c5} shows two drawings of $C_5\times C_2$, the second one making its symmetries more apparent.

\begin{theorem}\label{c2xc5thm} The graph $C_5 \times C_2$ has Gram dimension 4.
\end{theorem}

Throughout this section we set $G= C_5\times C_2=(V=[10],E)$. Clearly, $\gd(G)\ge 4$ since $K_4$ is a minor of $G$.
In order to show that $\gd(G)\le 4$,  we must show that $\gd(G,a)\le 4$ for any $a\in \SSS_{++}(G)$. Moreover, 
  in view of Corollary~\ref{lemgeneric}, it suffices to show this for all $a\in \SSS_{++}(G)$  satisfying the following `genericity' property:  For any Gram realization $\pb$ of $(G,a)$, 
  \begin{equation}\label{eqgene}
\rank \la \pb_C\ra \ge 3 \ \text{  for any circuit } C \text{  in } G. 
\end{equation}
From now on, we fix  $a\in \SSS_{++}(G)$ satisfying this genericity property.
Our objective  is to show that there exists a Gram realization of $(G,a)$ in $\oR^4$.

Again we use Theorem \ref{thmpsdstress} to construct an initial Gram realization $\pb$ of $(G,a)$. 
As stretched edge $e_0$, we choose the pair $e_0=(3,8)$ and we denote by $\hG=([8], \hE=E\cup\{(3,8)\})$ the extended graph obtained by adding the stretched pair $(3,8)$ to $G$.
By Theorem \ref{thmpsdstress}, we also have a stress matrix $\Om$ so that $(\hG,\pb,\Om)$ is a psd stressed framework.
Our objective is now  to  construct from $\pb$  another   Gram  realization $\bold{q}$ of $(G,a)$ lying in $\oR^4$.

\subsection{Additional useful  lemmas}\label{uaddlemmas}

First we deal with the case when $\rank \la \p_i,\p_j\ra=1$ for some pair $(i,j)$ of distinct nodes.
As $a\in \SSS_{++}(G)$, this can only happen when $(i,j)\not\in E$. 

\begin{lemma}\label{lem38}
If $\rank \la \p_i,\p_j\ra=1$ for some pair $(i,j)\not\in E$, then  there is a configuration in $\oR^4$ equivalent to $(G,\bold{p})$.
\end{lemma}

\begin{proof}
By assumption,  $\p_i=\epsilon \p_j$ for some scalar $\epsilon\ne 0$. 
Up to symmetry there are two cases to consider: (i) $(i,j)=(1,5)$ (two nodes at distance 2 in $G$), or (ii) $(i,j)=(1,6)$ (two nodes at distance 3). Consider first case (i) when $(i,j)=(1,5)$, so $\p_1=\epsilon \p_5$.
Set $V'=V\setminus \{1\}$. 
Let $G'=(V',E')$ be the graph on $V'$ obtained from $G$ by deleting node 1 and adding the edges $(2,5)$ and $ (5,9)$ (in other words, get $G'$ by identifying nodes 1 and 5 in $G$).
Let $X'$ be the Gram matrix of the vectors $\p_i$ ($i\in V'$) and define $a'=(X'_{jk})_{jk\in V'\cup E'} \in \SSS_+(G')$. First we show that $(G',a')$ has a Gram realization in $\oR^4$. For this, consider the graph $H$ obtained from $G$  by deleting both nodes 1 and 5. Then $G'$ is a subgraph of $\nabla H$ and thus
$\gd(G')\le \gd(\nabla H)=\gd(H)+1$. As $\tw(H)\le 2$ it follows that  $\gd(H)\le 3$ and thus $\gd(G')\le 4$.
Finally,   if ${\bf q}_{V'}$ is a Gram realization in $\oR^4$ of $(G',a')$ then, setting $\q_1=\epsilon\q_5$, we obtain a Gram realization ${\bf q}$ of $(G,a)$ in $\oR^4$.

Case (ii) is analogous, based on the fact that  the graph $H$ obtained from $G$ by deleting nodes 1 and 6 is a partial 2-tree.
\qed\end{proof}

\ignore{
\begin{corollary}\label{cor1node}
If there is a 1-node in the stressed graph $\SSS(\Om)$, then  there is a configuration in $\oR^4$ equivalent to $(G,p)$.
\end{corollary}

\begin{proof}
If $i$ is a 1-node then, by Lemma \ref{lem1node}, $\rank \la \p_i,\p_j\ra=1$ for some edge $ij$ of $\hG$. As $a\in \SSS_{++}(G)$, we must have $(i,j)=(3,8)$ and then we conclude using Lemma \ref{lem38}.
\qed\end{proof}
}

We  now  consider the case when the stressed graph might have a circuit as a connected component.

\begin{lemma}\label{lemnocircuit}
Let  $C$ be a circuit in $\hG$. If $C$ is   a connected component of  $\SSS(\Om)$, then $\rank \la \pb_C\ra \le 2$.
\end{lemma}

\begin{proof}
Directy, using Lemma \ref{contract2node} combined with Lemma \ref{base}.
\qed\end{proof}

Therefore, in view of the genericity assumption (\ref{eqgene}),  if a circuit $C$ is a connected component of the stressed graph, then $C$ cannot be a circuit in $G$ and thus  $C$ must contain the stretched pair $e_0=(3,8)$. 
The next result is useful to handle this case, treated in Corollary \ref{corcircuit38} below.

\begin{lemma}\label{lemcut}
Let $N_2(i)$ be the set of nodes at distance 2 from a given node $i$ in $G$.
If $\rank \la \pb_{N_2(i)}\ra \le 3$, then there is a configuration equivalent to $(G,\bold{p})$ in $\oR^4$.
\end{lemma}

\begin{proof}
Say, $i=1$ so that $N_2(1)= \{4,5,7,10\}$, cf. Figure~\ref{node1}.
Consider the set $S=\{2,3,6,9\}$ which is stable in $G$. 
Let $H$ denote the graph obtained from $G$ in the following way: For each node $i\in S$,  delete $i$ and add the clique on 
$N(i)$. One can verify that $H$ is contained in the clique 4-sum of the two  cliques $H_1$ and $H_2$ on the node sets 
$V_1=\{1,4,5,7,10\}$ and $V_2=\{4,5,7,8,10\}$, respectively.
By assumption, $\rank \la \pb_{V_1}\ra \le 4$ and $\rank \la \pb_{V_2}\ra \le 4$.
Therefore, one can apply an orthogonal transformation and find vectors $q_i\in \oR^4$ ($i\in V_1\cup V_2$) such that 
$\pb_{V_r}$ and $\bold{\q}_{V_r}$ have the same Gram matrix, for $r=1,2$.
Finally, as $V_1\cup V_2=V\setminus S$ and  the set $S$ is stable in $G$, one can  extend to a configuration ${\bf q}_V$ equivalent to $\pb_V$ by applying Lemma \ref{folding}.
\qed\end{proof} 

\begin{corollary}\label{corcircuit38}
If there is a circuit $C$   in $\hG$ containing the (stretched) edge $(3,8)$ such that $\rank \la \pb_C\ra  \le 2$, then  there is a configuration equivalent to $(G,\bold{p})$ in $\oR^4$.
\end{corollary}

\begin{proof}
If $|C|\ge 7$, pick $i\in V\setminus C$ and note that $\rank\la \pb_{-i}\ra  \le 4$. 
If $|C|=6$, pick a subset $S\subseteq V\setminus C$ of cardinality 2 which is stable in $G$, so that 
$\rank \la \pb_{V\setminus S}\ra \le 4$. In both cases we can conclude using Lemma~\ref{folding}.
Assume now that $|C|= 4$ or 5. 
In view of Lemma \ref{lemcut}, it suffices to check that  there exists a node $i$ for which 
$|C\cap N_2(i)| =3$.
For instance, for $C=(3,8,7,5)$, this holds  for node $i= 9$, and for $C=(3,8,10,9,1)$ this holds for $i=2$. Then,  Lemma~\ref{lemcut}  implies that $\rank \la \pb_{N_2(i)}\ra \le 3$. \qed\end{proof}

\begin{remark}\label{remcontract2node}
From now on, we  will assume that $\rank \la \p_i,\p_j\ra=2$ for all $i\ne j\in V$ (by Lemma \ref{lem38}).
 Hence there is no 1-node in the stressed graph.
Moreover, we will assume that no circuit $C$ of $\hG$ satisfies $\rank \la \bold{p}_C\ra \le 2$.
Therefore,    the stressed graph does not have a connected component which is a circuit (by  (\ref{eqgene}), Lemma \ref{lemnocircuit}  and Corollary \ref{corcircuit38}).
Hence we are guaranteed that after  contracting several 2-nodes   we do obtain a stressed framework (i.e, with a nonzero stress matrix).
\end{remark}

The next two lemmas settle the case when there are 
sufficiently many 2-nodes.

\begin{lemma}\label{lem42nodes}
If  there are at least four 2-nodes in the stressed graph $\SSS(\Om)$, then    there is a configuration equivalent to $(G,\pb)$ in $\oR^4$.
\end{lemma}

\begin{proof}
Let $I$ be a set of four 2-nodes in $\SSS(\Om)$. 
Hence, $\pb_I\subseteq \la \pb_{V\setminus I}\ra$ and thus it suffices to show that 
$\rank \la \pb_{V\setminus I}\ra \le 4$.

After contracting each of the four 2-nodes of $I$, we obtain a psd stressed framework 
$(\hG/ I, \pb_{V\setminus I }, \Om')$.
Indeed, we can apply Lemma \ref{contract2node} and obtain a nonzero psd stress matrix $\Om'$ in the contracted graph (recall Remark \ref{remcontract2node}).
If the support graph of $\Om'$ is not a clique,  Lemma \ref{base} implies  that
$\rank \la \pb_{V\setminus I}\ra  \le |V\setminus I|-2 = 4$.  

Assume now that $\SSS(\Om')$ is a clique on $T\subseteq V\setminus I$. Then  $\rank \la \pb_T\ra \le t-1$,
$|V\setminus (I\cup T)|=6-t$,  and $t=|T| \in \{3,4,5\}$. Indeed one cannot have $t=6$ since, after contracting the four 2-nodes, at least 4 edges are lost so that there remains at most $16- 4=12<15$ edges.
It suffices now  to  show that we can partition $V\setminus (I\cup T)$ as $S\cup S'$, where $S$ is stable in $G$ and 
$|S'|+t-1 \le 4$. Indeed, we then have $\rank \la \pb_{V\setminus S} \ra= \rank \la \pb_{T\cup S'}\ra \le  t-1+|S'|\le 4$ and we can conclude using Lemma \ref{folding}.
If $t=5$, then $|V\setminus (I\cup T)|=1$  and choose $S'=\emptyset$.
If $t=4$, then choose $S'\subseteq V\setminus (I\cup T)$ of cardinality 1. 
If $t=3$, then one can choose a stable set of cardinality 2 in $V\setminus (I\cup T)$ and $|S'|=1$.
\qed\end{proof}
 
\begin{lemma} \label{lem32nodes}
If there is at least one 0-node and at least  three 2-nodes in the stressed graph $\SSS(\Om)$,  then  there is a configuration equivalent to $(G,\pb)$ in $\oR^4$.
\end{lemma}

\begin{proof}
For $r=0,2$, let $V_r$ denote the set of $r$-nodes and set $n_r=|V_r|$. 
By assumption, $n_0\ge 1$ and we can assume $n_2=3$ (else apply Lemma \ref{lem42nodes}).
Set $W=V\setminus (V_0\cup V_2)$. 
After contracting the three 2-nodes in the stressed framework $(\hG,\pb,\Om)$, we get a stressed framework $(H,\pb_W,\Om')$ on $|W|=7-n_0$ nodes. 
Hence $n_0\le 4$ and  
 $\pb_{V_2}\subseteq  \la \pb_W\ra$.

Assume first that $\SSS(\Om')$ is not a clique. Then  $\rank \la \pb_W\ra \le |W|-2 =5-n_0$ by Lemma \ref{base}.
Now we can conclude using Lemma \ref{folding} since  in each of the cases: $n_0=1,2,3,4$, one can find a stable set $S\subseteq V_0$ such that $\rank \la \pb_{W\cup (V_0\setminus S)}\ra  \le 4$.

Assume now that $\SSS(\Om')$ is a clique. Then 
$\rank \la \pb_W\ra \le |W|-1=6-n_0$ by Lemma \ref{base}.
Note first that $n_0\ne 1,2$. Indeed, 
if $n_0=1$ then, after deleting the 0-node and  contracting the three 2-nodes,  
we have lost at least $3+3=6$ edges. Hence there remains at most $16-6=10$ edges in the stressed graph $\SSS(\Om')$, which therefore cannot be a clique on six nodes.
If $n_0=2$ then, after deleting the two 0-nodes and  contracting the three 2-nodes, we have lost at least $5 + 3=8$ edges. Hence there remain at most $16-8=8$ edges in the stressed graph $\SSS(\Om')$, which therefore cannot be a clique on five nodes. 
In each of the two cases $n_0=3,4$, one can find a stable set $S\subseteq V_0$ of cardinality 2 and thus 
$\rank \la \pb_{W\cup (V_0\setminus S)}\ra \le (6-n_0)+(n_0-2)=4$. Again conclude using Lemma \ref{folding}.
\qed\end{proof}

\subsection{Sketch of the proof}

In the proof we distinguish two cases: (i) when there is no 0-node, and (ii) when there is at least one 0-node, which are considered, respectively, in Sections \ref{secno0} and \ref{secyes0}.
In both cases the tools developed in the preceding section permit to conclude, except in one exceptional situation, occurring in case (ii). This execptional situation is when nodes 1,2,9 and 10 are 0-nodes and all edges of $\hG\setminus\{1,2,9,10\}$ are stressed. This situation needs a specific treatment which is done in Section \ref{secexceptionalcase}.

\subsection{There is no 0-node in the stressed graph}\label{secno0}
In this section we consider the case when  each node is stressed in  $\SSS(\Om)$, i.e.,  $w_{ii}\not=0$ for all $i\in [n]$.

\begin{lemma}\label{cyclelem} 
Assume that  all vertices are stressed in the stressed graph $\SSS(\Om)$ and that there exists a circuit  $C$ of length 4  in $G$ such that all edges in the cut $\delta(C)$ are stressed, i.e., $w_{ij}\ne 0$ for all edges $ij\in \hE$ with $i\in C$ and  $j\in V\setminus C$.  Then $\dim\langle\pb_V \rangle\le 4$.
\end{lemma}

\begin{proof}
 Up to symmetry, there are three types of  circuits $C$ of length 4  to consider:  (i) $C$ does not meet $\{3,8\}$, i.e., $C=(1,2,10,9)$;  or (ii) $C$ contains one of the two nodes 3,8, say node 8, and  it contains  a node adjacent  to the other one, i.e., node 3, like $C=(5,6,8,7)$;
 or (iii) $C$ contains one of 3,8 but has no node adjacent to the other one, like $C=(7,8,10,9)$. 

Consider first the case (i), when  $C=(1,2,10,9)$.  We show that the set $\pb_C$ spans $\pb_V$.
Using the equilibrium conditions at the nodes 1,2,9,10, we find that  $\p_3,\p_4,\p_7,\p_8 \in \langle \pb_C \ra$.
As 6 is not a 0-node, $w_{6i}\ne 0$ for some $i\in \{4,8\}$. Then, the equilibrium condition at node $i$ implies that 
$\p_6\in \la \pb_C\ra$.
Analogously for node 5.
\ignore{
If $w_{57}\not=0$, then the eq.~at 7 shows that $\p_5 \in \langle\p_1,\p_2,\p_9,\p_{10}\rangle$ and the eq.~condition at 6 that  $\p_6 \in \langle\p_1,\p_2,\p_9,\p_{10}\rangle$.

Suppose now that $w_{57}=0$. In this case, since there are no zero-nodes, 5 will be a 2-node and we can contract it. Then, one of $w_{46},w_{68}$ has to be nonzero. If $w_{46}\not=0$ then eq.~at 4 shows that $\p_6 \in \langle\p_1,\p_2,\p_9,\p_{10}\rangle$ and similarly for the other case. 
}

Case (ii) when $C=(5,6,8,7)$ can be treated in analogous manner. Just note that the equilibrium conditions applied to nodes 7,5,6 and 8 respectively, imply that $\p_9,\p_3,\p_4,\p_{10}\in \la \pb_C\ra$.   

We now consider case (iii) when $C=(7,8,10,9)$. Then one sees directly that $\p_1,\p_2,\p_5 \in \la \pb_C\ra$.
If $w_{24}\ne 0$, then the equilibrium conditions at nodes 2,3,6 imply that 
$\p_4,\p_3,\p_6 \in \la\pb_C\ra$ and thus $\la\pb_C\ra=\la \pb_V\ra $.  
Assume now that $w_{24}=0$, which implies $w_{34},w_{46}\ne 0$.
If $w_{13}\ne 0$,  then the equilibrium conditions at nodes 1,3,4 imply that  $\p_C$ spans $\p_3,\p_4,\p_6$ and we are done.
Assume now that $w_{24}=w_{13}=0$, so that  1,2,4 are 2-nodes. If there is one more 2-node then we are done by Lemma \ref{lem42nodes}.
Hence we can now assume that  $w_{ij}\ne 0$ whenever $(i,j)\ne (2,4)$ or $(1,3)$. 
After contracting the three 2-nodes 1,2,4 in the psd stressed framework $(\hG, \pb, \Om)$, we obtain a new psd stressed framework  on $V\setminus \{1,2,4\}$ where nodes 9, 10 have again degree 2. So contract these two nodes and get another psd  stressed framework on $V\setminus \{1,2,4,9,10\}$. Finally this implies 
$\rank \la \pb_V\ra =\rank \la \pb_{V\setminus \{1,2,4,9,10\}}\ra \le 4.$
\qed
\ignore{
In the first case, the eq.~conditions at 5,6,7,8 imply that $\p_3,\p_4,\p_9,\p_{10} \in \langle\p_5,\p_6,\p_7,\p_8\rangle$. If $w_{13}\not=0$ then eq.~at 3 gives that $\p_1 \in \langle\p_5,\p_6,\p_7,\p_8\rangle$ and then eq.~at 2 shows that $\p_2 \in \langle\p_5,\p_6,\p_7,\p_8\rangle$. On the other hand if $w_{13}=0$, since there are no zero-nodes, 1 is a 2-node and we can contract it.  Then, one of $w_{24},w_{2,10}$ is nonzero. If $w_{24}\not=0$ then eq.~at 4 shows that $\p_2 \in \langle\p_5,\p_6,\p_7,\p_8\rangle$ and similarly for the other case.

Assume now that $C=(7,8,9,10)$. In this case, doing a more tedious analysis we see that   $\langle\p\rangle=\langle\p_5,\p_6,\p_7,\p_8\rangle$.\qed\footnote{TO DO}
}
\end{proof}

In view of Lemma \ref{cyclelem}, we can now assume that, for each circuit $C$ of length 4 in $G$, there is at least one edge  $ij\in \delta (C)$ which is not stressed, i.e., $w_{ij}=0$.
It suffices now to show  that this implies the existence of at least four 2-nodes, as we can then conclude using Lemma \ref{lem42nodes}.

\medskip
For this let us enumerate the cuts $\delta (C)$ of the 4-circuits $C$  in $G$:\\
$\bullet$ For $C=(1,2,10,9)$, $\delta(C)=\{(1,3),(2,4),(7,9),(8,10)\}$.\\
$\bullet$ For $C=(7,9,10,8)$, $\delta(C)=\{(1,9), (2,10), (5,7), (6,8)\}.$\\
$\bullet $ For $C=(5,6,8,7)$, $\delta(C)=\{(7,9), (8,10), (3,5), (4,6)\}$.\\
$\bullet$ For $C=(3,5,6,4)$, $\delta(C)=\{(1,3), (2,4), (5,7),(6,8)\}$.\\
$\bullet$ For $C=(1,3,4,2)$, $\delta(C)=\{(3,5), (4,6), (1,9), (2, 10)\}$.\\
For instance,  $w_{24}=0$ implies that both 2 and 4 are 2-nodes, while $w_{13}=0$ implies that 1 is a 2-node.
One can easily verify that  there are at least four 2-nodes in $\SSS(\Om)$.

\subsection{There is at least one 0-node in the stressed graph}\label{secyes0}

Note that the mapping $\sigma: V\rightarrow V$ that permutes each of the pairs $(1,10)$, $(4,7)$, $(5,6)$, $(2,9)$ and $(3,8)$ is an automorphism of $G$. This can be easily seen  using the second drawing of $C_5 \times C_2$  in Figure~\ref{c2c5}. Hence, as nodes 3 and 8 are not 0-nodes, up to symmetry, it suffices to consider the following three cases:\\
$\bullet$ Node $1$ is a 0-node.\\
$\bullet$ Nodes 1, 10 are not 0-nodes and node 4 is a 0-node.\\
$\bullet$ Nodes 1, 10, 4, 7 are not 0-nodes and one of  5 or 2 is a 0-node.

\subsubsection{Node 1 is a 0-node.}\label{sec11}
\begin{figure}[h]
\centering
\includegraphics[scale=0.5]{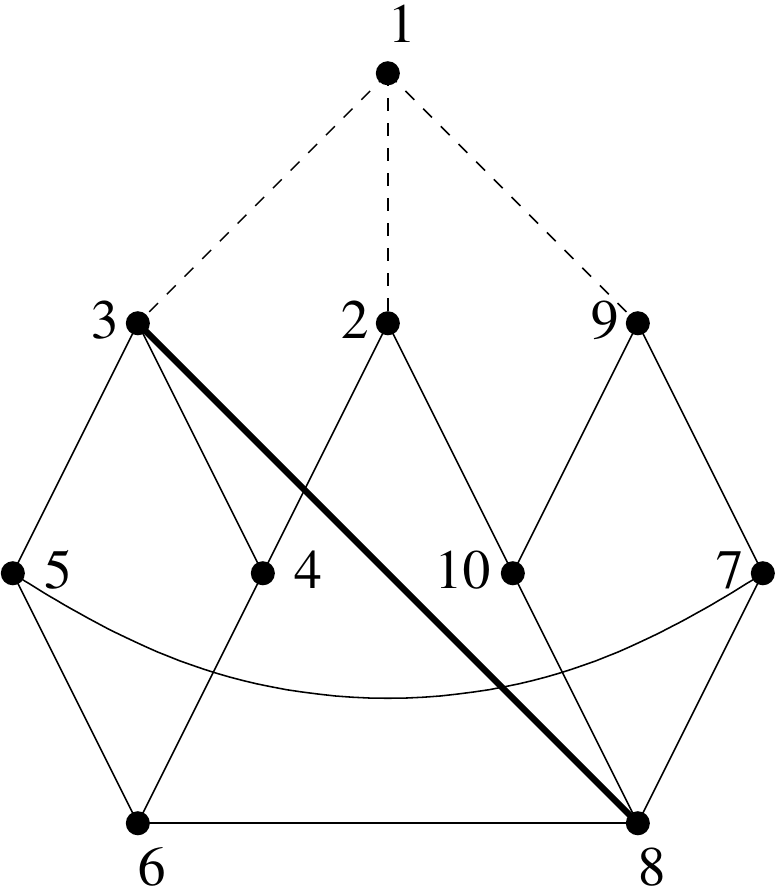}
\caption{A drawing of $\widehat{C_5 \times C_2}$ with 1 as the root node.}
\label{node1}
\end{figure}

It will be useful to use the drawing of $\hG$ from   Figure~\ref{node1}. There, the thick edge (3,8) is known to be stressed, the dotted edges are known to be non-stressed (i.e., $w_{ij}=0$),  while the other edges could be stressed or not.
In view of Lemma \ref{lem32nodes}, we can assume that there are at most two 2-nodes (else we are done).

\smallskip
Assume first that both nodes 2 and 9 are 0-nodes. Then node 10 too is a 0-node and each of nodes 4 and 7 is a 0- or 2-node.
If both 4,7 are 2-nodes, then all edges in the graph $G\backslash \{1,2,9,10\}$ are stressed. Hence we are in the {\em exceptional case,}  which we will consider in Section \ref{secexceptionalcase} below.
If 4 is a 0-node and 7 is a 2-node, then 3,7 must be the only 2-nodes and thus 6 is a 0-node. Hence,  the stressed graph is the circuit $C=(3,8,5,7)$, which implies $\rank\la {\bf p}_C\ra \le 2$ and thus we can conclude using Corollary \ref{corcircuit38}.
If 4 is a 2-node and 7 is a 0-node, then we find at least two more 2-nodes.
Finally, if both 4,7 are 0-nodes, then the stressed graph is the circuit $C=(3,8,6,5)$ and thus we can again conclude using 
Corollary \ref{corcircuit38}.

We can now assume that at least one of the two nodes 2,9 is a 2-node.
Then, node 3 has degree 3 in the stressed graph. (Indeed, if  3 is a 2-node, then 10 must be a 0-node (else we have three 2-nodes), which implies that 2,9 are 0-nodes, a contradiction.)
If exactly one of nodes 2,9 is stressed, one can easily see that there should be at least three 2-nodes. 
Finally consider the case when  both nodes 2,9 are stressed. Then they are the only 2-nodes which implies that all edges of $G\backslash 1$ are stressed. 
Set $I=\{4,5,8\}$. We show that   $\pb_I$ spans $\pb_{V\setminus \{1\}}$, so that $\pb_{\{1,4,5,8\}}$ spans $\pb_V$.
Indeed, the equilibrium conditions at 3 and 6 imply that $\p_3,\p_6\in \la \pb_I\ra$. 
Next, the equilibrium conditions at $4,5,2,9$ imply, respectively, that $\p_2\in \la \p_3,\p_4,\p_6\ra\subseteq \la\pb_I\ra$,
$\p_7\in \la \p_3,\p_5,\p_6\ra\subseteq\la\pb_I\ra$,
$\p_{10}\in\la \p_2,\p_4\ra\subseteq \la\pb_I\ra$, and
$\p_9\in \la \p_7,\p_{10}\ra\subseteq \la \pb_I\ra$.
This concludes the proof.

\subsubsection{Nodes 1, 10 are not 0-nodes and node 4 is a 0-node.}

\begin{figure}[h]
\centering
\includegraphics[scale=0.5]{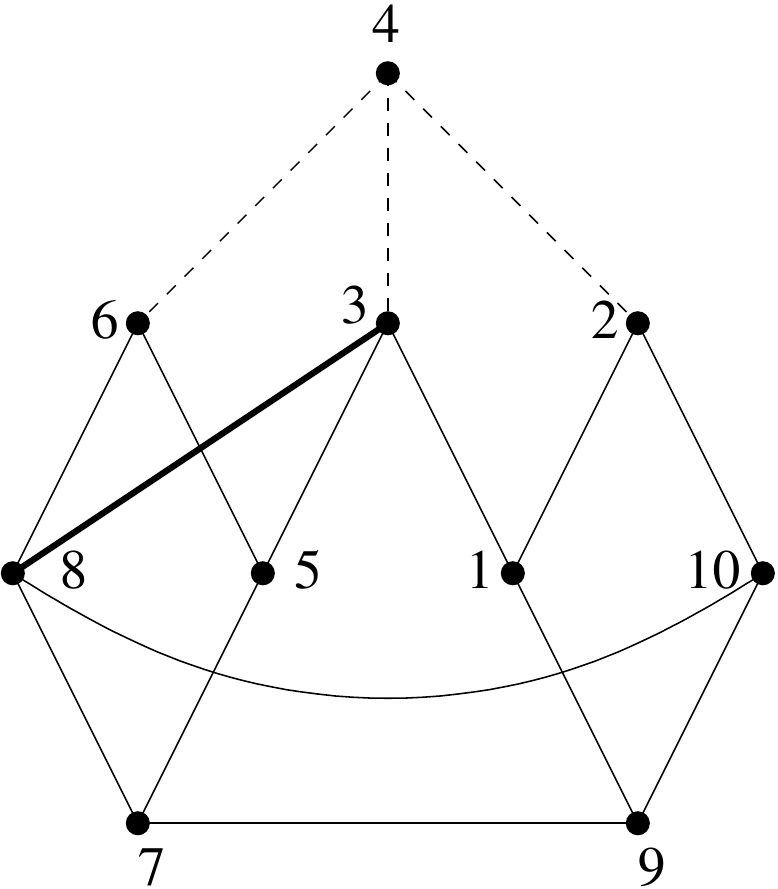}
\caption{A drawing of $\widehat{C_5 \times C_2}$ with 2 as the root node.}
\label{node2}
\end{figure}

It will be useful to use the drawing of $\hG$ from   Figure~\ref{node2}. We can assume that node 2 is a 2-node and that  node 3 has degree 3 in the stressed graph, since otherwise one would find at least  three 2-nodes. 
Consider first the case when 6 is a 2-node.

Then nodes 2 and 6 are the only 2-nodes which implies that all edges in the graph 
$G\backslash 4$ are stressed. Set $I=\{3,5,7,10\}$. We show that $\p_I$ spans $\p_{V\setminus \{4\}}$, and then we can conclude using Lemma \ref{folding}.
Indeed, the equilibrium conditions applied, respectively, to nodes
5,6,3,1,2 imply that 
$\p_6\in \la \pb_I\ra$,
$\p_8\in \la \p_5,\p_6\ra \subseteq \la \pb_I\ra$,
$\p_1\in \la \p_3,\p_5,\p_8\ra\subseteq \la \pb_I\ra$,
$\p_9\in \la \p_1,\p_7,\p_{10}\ra\subseteq \la\pb_I\ra$,
$\p_2\in \la \p_1,\p_{10}\ra\subseteq \la\pb_I\ra$.

Consider now the case when 6 is a 0-node.
Then 2 and 5 are the only 2-nodes so that all edges in the graph $G\backslash \{4,6\}$ are stressed.
Set $I=\{3,7,10\}$. We show that $\pb_I$ spans $\pb_{V\setminus \{4,6\}}$, and then we can again conclude using Lemma \ref{folding}.
Indeed the equilibrium conditions applied, respectively, at nodes 5,8,3,2,1 imply that
$\p_5, \p_8 \in  \la \pb_I\ra$,
$\p_1\in \la \p_3,\p_5,\p_8\ra\subseteq \la \pb_I\ra$,
$\p_2\in \la \p_1,\p_{10}\ra\subseteq \la \pb_I\ra$,
$\p_9\in \la \p_2,\p_8,\p_{10}\ra\subseteq \la\pb_I\ra$.

\subsubsection{Nodes 1, 4, 7, 10  are not 0-nodes and node 5 or 2  is a 0-node.}

\begin{figure}[h]
\centering
\includegraphics[scale=0.5]{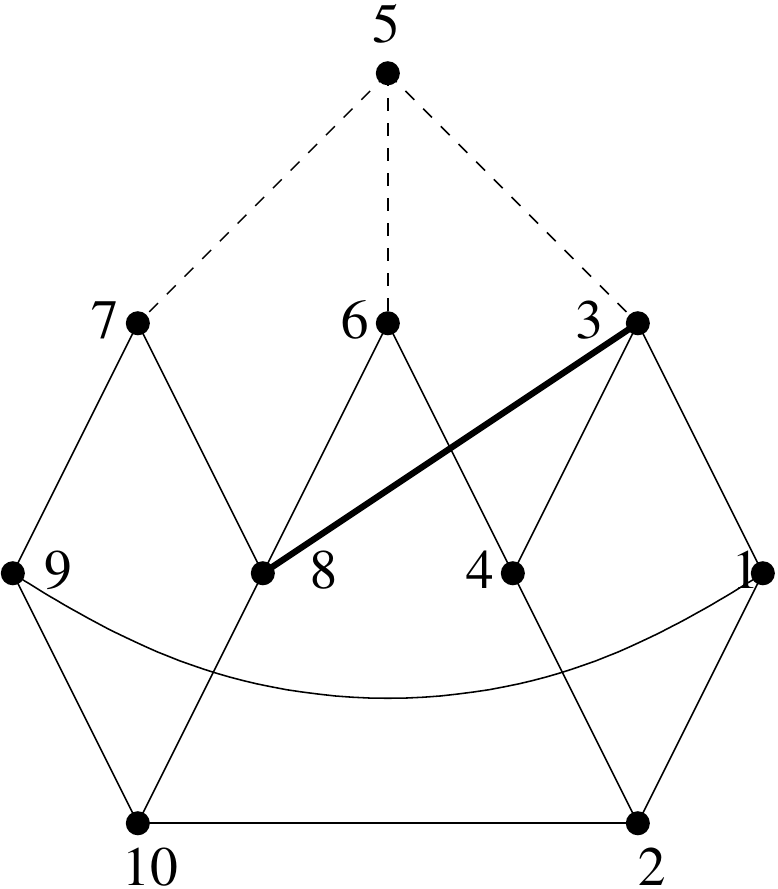}
\caption{A drawing of $\widehat{C_5 \times C_2}$ with 3 as the root node.}
\label{node3}
\end{figure}

It will be useful to use the drawing of $\hG$ from   Figure~\ref{node3}. We assume that nodes 1,4,7,10 are not 0-nodes. Consider first the case when node 5 is a 0-node.
Then node 7 is a 2-node. 

If  node 6 is a 2-node, then 6 and 7 are the only 2-nodes and thus all edges of the graph $G\backslash 5$ are stressed. Setting $I=\{1,2,4,8\}$, one can verify that $\pb_I$ spans $\pb_{V\setminus \{5\}}$ and then one can conclude using Lemma \ref{folding}.

If  node 6 is a 0-node, then
 nodes 4 and 7 are the only 2-nodes and thus all edges in the graph $G\backslash \{5,6\}$ are stressed. Setting $I=\{2,3,9\}$, one can verify that $\pb_I$ spans $\pb_{V\setminus \{5,6\}}$. Thus $\pb_{\{2,3,9,6\}}$ spans $\pb_{V\setminus \{5\}}$ and one can again conclude using Lemma \ref{folding}.

Consider finally the case when nodes 1,4,7,10, 5 and 6 are not 0-nodes and node 2 is a 0-node.
As node 2 is adjacent to nodes 1, 4 and 10 in $G$, we find three 2-nodes and thus we are done.

\subsection{The exceptional case}\label{secexceptionalcase}

In this section we consider the following case which was left open in the first case considered in Section \ref{sec11}:
 Nodes 1, 2, 9 and 10 are 0-nodes and all edges of the graph $\hG\backslash \{1,2,9,10\}$ are stressed. 

Then, nodes 4 and 7 are 2-nodes in the stressed graph. After contracting both nodes 4,7, we obtain a stressed graph which is the complete graph on 4 nodes. Hence, using Lemma \ref{base},
 we can conclude that  $\dim \la\pb_{V_1}\ra\le 3$, where  $V_1=  V\setminus  \{1,2,9,10\}$.
  Among the nodes 1, 2, 9 and 10, we can  find a stable set of size 2. Hence, if $\dim\la\pb_{V_1}\ra\le 2$ then, 
    using Lemma \ref{folding},  we can  find an equivalent configuration in dimension $2+2 =4$ and we are done.
From now on we assume that   
\begin{equation}\label{dimstressed} \dim\la \pb_{V_1}\ra=3.
\end{equation}
In this case it is not clear how to fold $\pb$ in $\oR^4$.
 In order to settle this case, we proceed as in Belk \cite{Belk}: We fix (or {\em pin}) the vectors $\p_i$ labeling the nodes $i\in V_1$ and we search for another set of vectors  $\p'_i$  labeling the nodes $i\in V_2=V\setminus V_1=\{1,2,9,10\}$ so that $\pb_{V_1}\cup \pb'_{V_2}$ can be folded into  $\oR^4$.
 Again, our starting point is to get  such new  vectors $\p'_i$ ($i\in V_2$) which, together with $\pb_{V_1}$, provide a Gram realization of $(G,a)$, by stretching along a second pair $e'$; namely we stretch the pair $e'=(4,9)\in V_1\times V_2$.
As in So and Ye \cite{SY06}, this  configuration $\pb'_{V_2}$ is again obtained by solving a semidefinite program; details are given below. 
 
 \subsubsection{Computing $\pb'_{V_2}$  via semidefinite programming.}

Let  $E[V_2]$ denote the set of edges of $G$ contained in $V_2$ and   let $E[V_1,V_2]$ denote the set of edges $(i,j)\in E$ with $i\in V_1$, $j\in V_2$. 
Moreover, set $|V_1|=n_1 \ge |V_2|=n_2$,  so the configuration $\pb_{V_1}$ lies in $ \oR^{n_1}$. (Here $n_1=6$, $n_2=4$).  We now search for a new configuration $\pb'_{V_2}$ by  stretching along the pair $(4,9)$. 
For this we use  the following semidefinite program:
\begin{equation}\label{SDPPpin}
\begin{array}{lll}
\max\ \la F_{49},Z\ra \ \text{\rm  such that } &    \la F_{ij},Z\ra =a_{ij} &  \forall  ij \in E[V_1,V_2]\\
& \la E_{ij}, Z\ra =a_{ij}  & \forall ij \in V_2\cup E[V_2]\\
&  \la E_{ij},Z\ra =0 & \forall i< j, i,j\in V_1\\
& \la E_{ii},Z\ra =1& \forall i\in V_1\\
&Z\succeq 0.
\end{array}
\end{equation}
Here,  $E_{ij}=(e_ie_j^T+e_je_i^T)/2 \in \SSS^{n_1+n_2} $, where $e_i$ ($ i\in [n_1+ n_2]$) are the standard unit vectors in  $\oR^{n_1+n_2}$. 
Moreover, for $i\in V_1$, $j\in V_2$, $F_{ij}= (p_i'e_j^T+e_j(p_i')^T)/2$, after setting
$p'_i=(\p_i,0)\in \oR^{n_1+n_2}$.

\medskip
Consider now  a matrix $Z$ feasible for~(\ref{SDPPpin}). Then $Z$ can be written in the block form
$Z=\left(\begin{matrix} I_{n_1} & Y \cr Y^T & X\end{matrix}\right)$, and let $\y_i\in \oR^{n_1}$ ($i\in V_2$) denote the columns of $Y$. The condition $Z\succeq 0$ is equivalent to $X-Y^TY\succeq 0$. Say,  $X-Y^TY$ is the Gram matrix of 
the vectors $\z_i\in \oR^{n_2}$ ($i\in V_2$). 
For  $i\in V_2$,  set  $\p'_i=(\y_i,\z_i)\in \oR^{n_1+n_2}$. Then   $X$ is the Gram matrix of the vectors  $\p'_i  \ (i \in V_2)$. 

For $i\in V_1$, $j\in V_2$ we have that   $\la F_{ij},Z\ra =(\p_i,0)^T(\y_j,\z_j)=(\p'_i)^T \p'_j.$ 
Moreover, for $i,j\in V_2$,  we have that   $\la {E}_{ij},Z\ra =X_{ij}=(\p'_i)^T\p'_j$.

Therefore, the linear conditions $\la F_{ij},Z\ra=a_{ij}$ for $ij \in E[V_1,V_2] $     and $\la E_{ij},Z\ra =a_{ij} $ for $ij\in V_2\cup E[V_2]$ imply that the vectors 
$\p'_i$ ($i\in V_1\cup V_2$) 
 form a Gram realization of $(G,a)$.

\medskip
We now consider the  dual semidefinite program of (\ref{SDPPpin}) which, as we see in Lemma \ref{lemOmp} below,  will give us some  equilibrium conditions  on the new vectors $\p'_i$ ($i\in V_2$). 
The dual program  involves scalar variables $w'_{ij}$ (for $ij\in E[V_1,V_2]\cup V_2\cup E[V_2]$) and a matrix
$U'=\left(\begin{matrix} U & 0 \cr 0 & 0\end{matrix}\right)$, and  it reads:

\begin{equation}\label{SDPDpin}
\begin{array}{ll}
\min & \displaystyle  \la I_{n_1}, U\ra + \sum_{ij \in E[V_1,V_2]} w'_{ij} a_{ij} +\sum_{ij \in V_2\cup E[V_2]} w'_{ij} a_{ij} \\
\text{such that }& \displaystyle \Om' = -F_{49} + U' +\sum_{ij \in E[V_1,V_2] } w'_{ij}F_{ij} +\sum_{ij\in V_2\cup E[V_2]} w'_{ij} E_{ij}\succeq 0.
\end{array}
\end{equation}

Since the primal program~(\ref{SDPPpin}) is bounded and  the dual program (\ref{SDPDpin}) is strictly feasible it follows that  program (\ref{SDPPpin}) has an optimal solution
$Z$.
Let $\p'_i \in \oR^{n_1+n_2}$ ($i\in V_1\cup V_2$) be the vectors as defined above, which thus form a Gram realization of 
$(G,a)$.

\begin{lemma}\label{lemOmp}
There exists a nonzero matrix $\Om'=(w'_{ij}) \succeq 0$ satisfying the optimality condition $Z\Om'=0$ and 
 the following conditions on its support: 
 \begin{equation}\label{support}
 \begin{array}{l}
  w'_{ij}=0 \ \ \forall  (i,j)\in (V_1\times V_2) \setminus  (E[V_1,V_2] \cup \{(4,9)\}), \\
w'_{ij}=0  \ \ \forall i\ne j\in V_2,  (i,j) \not\in E[V_2].
\end{array}
\end{equation}
 Moreover, the following equilibrium conditions hold:
\begin{equation}\label{newequi}
 w'_{ii}\p'_i+\sum_{j\in V_1\cup V_2\mid ij\in  E\cup \{(4,9)\}  }  
 w'_{ij}\p'_j=0\ \ \forall i\in V_2
 \end{equation}
 and $w'_{ij}\ne 0$ for some $ij\in V_2\cup E[V_2]$.
 Furthermore,  a node $i\in V_2$ is a 0-node, i.e., $w'_{ij}=0$ for all $j\in V_1\cup V_2$,  if and only if $w'_{ii}=0$.
 \end{lemma}

\begin{proof}
If the primal program (\ref{SDPPpin}) is strictly feasible,  then (\ref{SDPDpin}) has an optimal solution $\Om'$ which satisfies $Z\Om'=0$ and  (\ref{support}) (with $w'_{49}=-1$).
Otherwise, if (\ref{SDPPpin}) is feasible but not strictly feasible then, using Farkas' lemma (Lemma \ref{lemfarkas}),
we again find a matrix $\Om'\succeq 0$ satisfying $Z\Om'=0$ and (\ref{support}) (now with $w'_{49}=0$).
We now indicate how to derive (\ref{newequi}) from the condition $Z\Om'=0$. 

For this write the matrices $Z$ and $\Om'$ in block form
$$Z=\left(\begin{matrix} I_{n_1} & Y\cr Y^T & X \end{matrix}\right),\ \ \Om'=\left(\begin{matrix} \Om'_{1}& \Om'_{12}\cr 
(\Om'_{12})^T & \Om'_2\end{matrix}\right).$$ From  $Z\Om'=0$, we derive 
$Y^T\Om'_{12}+X\Om'_2=0$ and  $\Om'_{12}+Y\Om'_2=0$. First this implies 
$(X-Y^TY)\Om'_2=0$ which in turn  implies that the $V_2$-coordinates of the vectors on  the left hand side in (\ref{newequi}) are equal to 0. 
Second,  the condition $\Om'_{12}+Y\Om'_2=0$ together with expressing  $\Om'_{12}= 
\sum_{ij\in E[V_1,V_2]\cup \{(4,9)\} } w'_{ij} \p'_i e_j^T $,  implies that the $V_1$-coordinates of the vectors on  the left hand side in (\ref{newequi}) are 0. Thus (\ref{newequi}) holds.
Finally, we verify that $\Om'_2\ne 0$. Indeed, $\Om'_2=0$ implies $\Om'_{12}=0$ and thus $\Om'=0$ since
$0=\la Z,\Om'\ra = \la I_{n_1}, \Om'_1\ra$.
\qed 
\end{proof}

\subsubsection{Folding $\pb'$ into $\oR^4$.}

 We now use the above configuration $\p'$  and the equilibrium conditions (\ref{newequi}) at the nodes of $V_2$  to construct  a Gram realization of $(G,a)$ in $\oR^4$.
By construction, $\p'_i=(\p_i,0)$ for $i\in V_1$. 
Note that  no node $i\in V_2$ is a 1-node with respect to the new stress $\Omega'$ (recall Lemma \ref{lem38}).
Let us point out again   that  Lemma~\ref{lemOmp}  does not guarantee that  $w'_{49}\ne 0$  (as opposed to relation (\ref{patternOm0}) in Theorem~\ref{thmpsdstress}).
 
By assumption nodes 1,2,9 and 10 are 0-nodes and all other edges of the graph $\hG\setminus \{1,2,9,10\}$ are stressed w.r.t. the old stress matrix $\Omega$.  We begin with the following easy   observation about $\pb'_{V_1}$. 
 
  \begin{lemma}\label{lempV1}
 $\dim \la \p'_4,\p'_7,\p'_8\ra = \dim \la \p'_3,\p'_4,\p'_8\ra =3$. 
 \end{lemma}
 
 \begin{proof}
It is easy to see that each of these sets  spans $\pb'_{V_1}$ and $\rank \la \pb'_{V_1}\ra =\rank \la \pb_{V_1}\ra=3$ by ~(\ref{dimstressed}).
 \qed\end{proof}

\begin{figure}[h] 
\centering \includegraphics[scale=0.4]{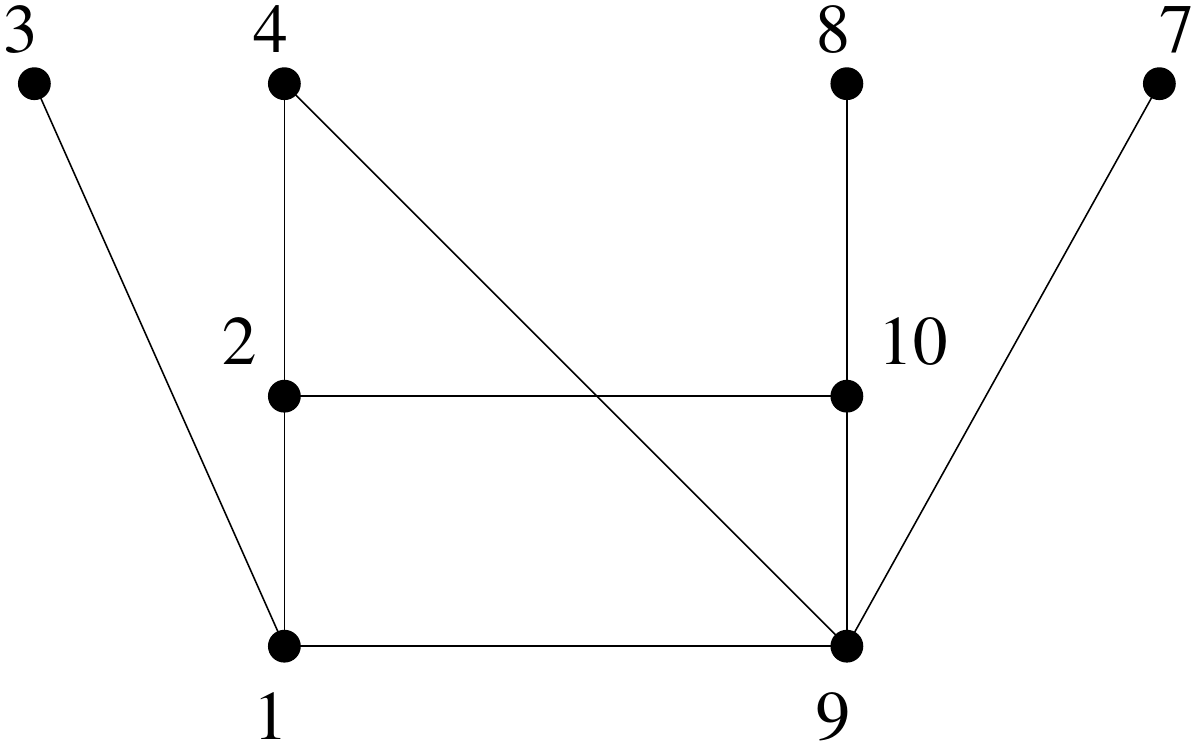}
\caption{The graph containing the support of the new stress matrix $\Om'$}
\label{exccase} 
\end{figure}

 As an immediate corollary we    may assume that 
\begin{equation}\label{as}
 \p'_i\not\in \la \pb'_{V_1}\ra \ \forall i\in V_2
\end{equation}
Indeed, if there exists $i\in V_2$ satisfying  $\p'_i \in \la \pb'_{V_1}\ra$ then we can  find a stable set of size two in $V_2\setminus \{i\}$ and   using Lemma \ref{folding} we can construct  an equivalent configuration in $\oR^4$.
Therefore,  we can assume that  at most two nodes in $V_2$ are 0-nodes in $\SSS(\Omega{'})$ since, by construction, for  the new stress matrix $\Omega'$  there exists $ij \in V_2\cup E[V_2]$ such that  $w'_{ij}\not=0$.  This guides our discussion below.   Figure~\ref{exccase}  shows the graph 
containing the support of the new stress matrix $\Omega'$.

\vspace*{-0.5cm}
\subsubsection*{There are two 0-nodes in $V_2$.} The cases when either 2,9, or 1,10, are 0-nodes are excluded (since then one would have a 1-node). 
 If nodes 1 and 9 are 0-nodes, then the equilibrium conditions at nodes 2 and 10 imply that
 $\p'_4,\p'_8\in \la \p'_2,\p'_{10}\ra$ and  by Lemma~\ref{lempV1} we have that  $\la \p'_2,\p'_{10}\ra =\la  \p'_4,\p'_8\ra \subseteq \la \pb'_{V_1}\ra$, contradicting (\ref{as}). The case when nodes 9,10 are 0-nodes is similar.

 Finally assume that nodes 1,2 are 0-nodes (the case when 2,10 are 0-nodes is analogous).
 As $w'_{8, 10}\ne 0$, the equilibrium condition at node 10 implies that $\p'_8\in \la \p'_9,\p'_{10}\ra$.
 If $w'_{49}=0$ then the equilibrium condition at node 9 implies that $\p'_7\in \la \p'_9,\p'_{10}\ra$.
 Hence $\la \p'_7,\p'_8\ra \subseteq  \la \p'_9,\p'_{10}\ra$, thus equality holds, contradicting (\ref{as}).
 If $w'_{49}\ne 0$, then $\p'_4\in \la \p'_7,\p'_9,\p'_{10}\ra$ and thus 
  $\la \p'_4, \p'_7,\p'_8\ra \subseteq \la \p'_7,\p'_9,\p'_{10}\ra$.
 Hence equality holds (by Lemma \ref{lempV1}), contradicting again (\ref{as}).
 
 \subsubsection*{There is one 0-node in $V_2$.}
 Suppose first 9 is the only 0-node in $V_2$. The equilibrium conditions at nodes 1 and 10 imply that
 $\p'_1\in \la \p'_3,\p'_2\ra$ and $\p'_{10}\in \la \p'_2, \p'_8\ra$. Hence $\la \p'_1,\p'_{10}\ra \subseteq \la \p'_{V_1}, \p'_2\ra$
 and thus $\rank \la \pb'_{V\setminus \{9\}}\ra = 4$. Then conclude using Lemma~\ref{folding}.
 
 Suppose now that node 1 is the only 0-node (the cases when 2 or 10 is 0-node are analogous).
 The equilibrium conditions at nodes 2 and 9 imply that 
 $\p'_2 \in \la \p'_4,\p'_{10}\ra$ and $\p'_9\in \la \p'_4,\p'_7,\p'_{10}\ra$.
 Hence, $\la \p'_2,\p'_9 \ra \subseteq \la \p'_{V_1},\p'_{10}\ra$ and we can conclude using Lemma \ref{folding}.
 
 \subsubsection*{There is no 0-node in $V_2$.}
 
 We can assume that  $w'_{ij}\ne 0$ for some $(i,j)\in V_1\times V_2$ for otherwise we get the stressed circuit $C=(1,2,10,9)$, thus with  $\rank \la \pb'_C\ra=2$, contradicting Corollary~\ref{lemgeneric}.
 We show that $\rank \la \pb'_V\ra=4$. 
 For this we discuss according to how many parameters are equal to zero  among $w'_{13},w'_{24}, w'_{8,10}$.
 If none is zero, then  the equilibrium conditions at nodes 1,2 and 10 imply that 
 $\p'_3,\p'_4, \p'_8\in \la \pb'_{V_2}\ra$ and thus Lemma~\ref{lempV1} implies that  
 $\dim ( \la \pb'_{V_1}\ra \cap \la \pb'_{V_2}\ra) \ge 3$. Therefore, 
  $\dim \la \pb'_{V_1}, \pb'_{V_2}\ra = \dim 
 \la \pb'_{V_1}\ra + \dim \la \pb'_{V_2}\ra - \dim (\la \pb'_{V_1}\ra \cap \la \pb'_{V_2}\ra)
  \le 3 + 4 -3 =4$.
  
  Assume now that  (say) $w'_{13}=0$, $w'_{24},w'_{8,10}\ne 0$.
  Then $\rank \la \pb'_{V_2}\ra \le 3$ (using the equilibrium condition at node 1).
  As $w'_{24},w'_{8,10}\ne 0$, we know that $\p'_4,\p'_8\in \la \pb'_{V_2}\ra$.
  Hence $  \dim (\la \pb'_{V_1}\ra \cap \la \pb'_{V_2}\ra)  \ge 2$ and thus 
  $\dim \la \pb'_{V_1}, \pb'_{V_2}\ra  
  \le 3 + 3 - 2=4$.
  
  Assume now (say) that $w'_{13}=w'_{24}=0$, $w'_{8,10}\ne 0$.
  Then the equilbrium conditions at nodes  1 and 2 imply that $\rank \la \pb'_{V_2}\ra\le 2$.
  Moreover, $\p'_8\in \la \pb'_{V_2}\ra$. Hence $  \dim (\la \pb'_{V_1}\ra \cap \la \pb'_{V_2}\ra) \ge 1$  
  and thus $\dim \la \pb'_{V_1}, \pb'_{V_2}\ra \le 3 + 2 -1 =4$.
  
    Finally assume now that $w'_{13}=w'_{24}=w'_{8,10}=0$. Then $\rank \la\pb'_{V_2}\ra= 2$.
   Moreover,  at least one of $w'_{49},w'_{79}$ is nonzero. Hence $  \dim (\la \pb'_{V_1}\ra \cap \la \pb'_{V_2}\ra)  \ge 1$ and thus
    $\dim \la \pb'_{V_1}, \pb_{V_2}\ra \le 3 + 2 -1 =4$.

 \section{Concluding remarks} 
 \label{sec:concluding}
        


One of the main   contributions of this paper is the proof that 
 $\gd(C_5 \times C_2)\le~4$, an inequality  which  underlies  the characterization of graphs with Gram dimension at most four. As already explained we obtain as corollaries the inequalities $\fedm(C_5\times C_2)\le 3$ of \cite{Belk} and $\nu^=(C_5\times C_2)\le 4$ of \cite{H03}.
 
 
 Although our  proof of  the inequality $\gd(C_5 \times C_2)\le 4$  goes roughly along  the same lines as the   proof of the inequality $\fedm(C_5 \times C_2) \le 3$ given in \cite{Belk}, there are important differences and we believe  that  our   proof 
     is simpler. 
 This is due in particular  to the fact that we introduce  a number of new  auxiliary lemmas (cf. Sections~\ref{ulemmas} and \ref{uaddlemmas}) that enable us to  deal more efficiently with  the case checking which constitutes the most tedious part of the proof.
Furthermore, the use of semidefinite programming  to construct a stress matrix permits to eliminate some case checking since, as was  already noted in~\cite{SY06} (in the context of Euclidean realizations), the stress is nonzero along the stretched pair of vertices. 
Additionally,  our analysis complements and at several occasions  even corrects the proof in~\cite{Belk}. As an example, the case when two vectors labeling two non-adjacent nodes are parallel in not discussed in~\cite{Belk}; this leads to some additional case checking which we address in Lemma~\ref{lem38}.


\medskip
As the class of graphs with $\gd(G)\le k$ is closed under taking minors, it follows from the general theory of Robertson and Seymour \cite{RS} that there exists a polynomial time algorithm for testing  $\gd(G)\le k$. 
For $k\le 4$ the forbidden minors  for $\gd(G)\le k$  are known, hence one can make this polynomial time algorithm explicit.
We refer to \cite[\S4.2.5]{So} for details on how to check $\gd(G)\le 4$ or, equivalently, $\fedm(G)\le 3$.

The next algorithmic question is how to construct a Gram representation in $\oR^4$ of a given partial matrix  $a\in \SSS_+(G)$ when $G$ has Gram dimension 4. As explained in \cite[\S4.2.4,\S4.2.5]{So}, the first step consists of finding  a graph $G'$ containing $G$ as a subgraph and such that $G'$ is a clique sum of copies of $V_8$, $C_5\times C_2$ and chordal graphs with tree-width at most 3. 
Then, if a  psd completion $A$ of $a$ is available, 
it suffices to deal with each of these components separately. Such a psd completion can be computed approximately by solving a semidefinite program.
Chordal components are easy to deal with in view of the general results on psd completions in the chordal case.
For the components of the form $V_8$ or $C_5\times C_2$ one has to go through the steps of the proof to get a new Gram representation in $\oR^4$.
The basic ingredient of  our proof is the {\em existence} of a primal-dual pair of optimal solutions to the programs~(\ref{SDPP}) and (\ref{SDPD}). Under appropriate genericity assumptions, the existence of such a pair of optimal solutions is guaranteed by standard results  of semidefinite programming duality theory (cf. Section~\ref{sec:sdpformulation}). Also in the case when the primal program is not strictly feasible,  we can still guarantee the existence of a psd stress matrix; our  proof  uses Farkas' lemma and  is  simpler  than the proof in~\cite{So} of the analogous result  in the context of Euclidean realizations. However, in the case of the graph $C_5\times C_2$, we must make an additional  genericity assumption on the vector  $a\in \SSS_{++}(G)$ (namely, that  the configuration restricted to any circuit is not coplanar). This is problematic since the folding procedure apparently breaks down for non-generic configurations; note that this issue also arises in the case of Euclidean embeddings since an analogous genericity assumption is made in \cite{Belk}, although it is not discussed  in the algorithmic approach of \cite{So,SY06}.
Moreover, the above procedure relies on solving several semidefinite programs, which however cannot be solved exactly in general, but only to some given precision. 
This thus excludes the possibility of turning  the proof  into an efficient algorithm for computing exact Gram representations in $\oR^4$.

\medskip
We conclude with the following question about the relation between the two parameters $\gd(G)$ and $\fedm(G)$, which has been  left open: Prove or disprove the inequality:
$$\fedm(\nabla G)\le \fedm(G)+1.$$
 
\medskip\noindent
{\bf Acknowledgements.} We thank M. E.-Nagy for useful discussions and  A. Schrijver for his suggestions for the  proof of Theorem  \ref{lemedmdim}.

\end{document}

